\documentclass[3p]{article}
\usepackage{graphicx}
\usepackage{xcolor}
\usepackage{latexsym,bm}
\usepackage{url}
\usepackage{xcolor}
\usepackage{rotating}
\usepackage{pdflscape}
\definecolor{newcolor}{rgb}{.8,.349,.1}
\usepackage{enumitem}

\usepackage{amsmath}	
\usepackage{subfig}
\usepackage{amssymb} 
\usepackage{amsfonts}
\usepackage{amsbsy}
\usepackage{mathtools}
\usepackage{url}
\usepackage{multirow}

\usepackage{caption}
\usepackage{varwidth}
\usepackage{chngcntr}
\usepackage{etoolbox}
\usepackage[toc,page]{appendix}
\usepackage{hyperref}
\hypersetup{
    colorlinks,
    linkcolor={black!50!black},
    citecolor={black!50!black},
    urlcolor={black!80!black}
}

\usepackage[colorinlistoftodos]{todonotes}
\usepackage{empheq}
\usepackage{indentfirst}
\usepackage[below]{placeins}
\usepackage{pdfpages}
\usepackage{comment}
\usepackage{booktabs}

\newcommand{\aposteriori}{\textit{a posteriori }}

\newcommand{\bu}{\boldsymbol{u}}
\newcommand{\bU}{\mathbf{U}}

\newcommand{\vel}{u}
\newcommand{\urel}{w}
\newcommand{\bvel}{\boldsymbol{u}}
\newcommand{\burel}{\boldsymbol{w}}

\newcommand{\bg}{\boldsymbol{g}}
\newcommand{\ff}{\mathbf{f}}
\renewcommand{\gg}{\mathbf{g}}
\newcommand{\hh}{\mathbf{h}}
\newcommand{\bQ}{\mathbf{Q}}
\newcommand{\bV}{\mathbf{V}}
\newcommand{\bF}{\mathbf{F}}
\newcommand{\bB}[1]{\mathbf{B}_{#1}}
\newcommand{\bS}{\mathbf{S}}
\newcommand{\bcalS}{\boldsymbol{\mathcal{S}}}
\newcommand{\bx}{\mathbf{x}}
\newcommand{\bq}{\mathbf{q}}
\newcommand{\bn}{\mathbf{n}}

\newcommand{\bA}{\mathbf{A}}
\newcommand{\bC}{\mathbf{C}}
\newcommand{\bR}{\mathbf{R}}

\newcommand{\bpsi}{\boldsymbol{\psi}}
\newcommand{\velr}{\vel^{r}}
\newcommand{\velt}{\vel^{\theta}}
\newcommand{\velz}{\vel^{z}}
\newcommand{\urelr}{\urel^{r}}
\newcommand{\urelt}{\urel^{\theta}}
\newcommand{\urelz}{\urel^{z}}

\newcommand{\dalpha}[1]{\der \alpha_{#1}}

\newcommand{\der}{\partial}
\newcommand{\dt}{\der t}

\newcommand{\dx}[1]{\der x_{#1}}
\newcommand{\ddx}[1]{\mathrm{d}{#1}}

\newcommand{\ds}{\displaystyle}

\allowdisplaybreaks

\newcommand{\footremember}[2]{%
	\footnote{#2}
	\newcounter{#1}
	\setcounter{#1}{\value{footnote}}%
}
\newcommand{\footrecall}[1]{%
	\footnotemark[\value{#1}]%
} 
\title{High-order ADER Discontinuous Galerkin schemes for a symmetric hyperbolic model of compressible barotropic two-fluid flows} 
\author{%
	Laura R\'io-Mart\'in\setcounter{footnote}{2}\footremember{4}{Department of Information Engineering and Computer Science, University of Trento, via Sommarive 9, Povo, 38123 Trento, Italy}\setcounter{footnote}{1}\footremember{2}{Laboratory of Applied Mathematics, DICAM, University of Trento, Via Mesiano 77, 38123 Trento, Italy}\setcounter{footnote}{3}\footnote{laura.delrio@unitn.it}%
	\and Michael Dumbser\footrecall{2} \setcounter{footnote}{0}\footnote{Corresponding author, michael.dumbser@unitn.it}
}
\date{}

\providecommand{\keywords}[1]
{
	\small	
	\textbf{\textit{Keywords---}} #1
}
	
\begin{document}
	\maketitle
\begin{abstract} 
This paper presents a high-order discontinuous Galerkin finite element method to solve the barotropic version of the conservative symmetric hyperbolic and thermodynamically compatible (SHTC) model of compressible two-phase flow, introduced by Romenski \textit{et al.} in~\cite{Romenski2007TwoPhase,Romenski2010TwoPhase}, 
in multiple space dimensions. In the absence of algebraic source terms, the model is endowed with a curl constraint on the relative velocity field. 
In this paper, the hyperbolicity of the system is studied for the first time in the multidimensional case, showing that the original model is only weakly hyperbolic in multiple space dimensions. To restore strong hyperbolicity, two different methodologies are used: i) the explicit symmetrization of the system, which can be achieved by adding terms that contain linear combinations of the curl involution, similar to the Godunov-Powell terms in the MHD equations; ii) the use of the hyperbolic generalized Lagrangian multiplier (GLM) curl-cleaning approach forwarded. The PDE system is solved using a high-order ADER discontinuous Galerkin method with \aposteriori sub-cell finite volume limiter to deal with shock waves and the steep gradients in the volume fraction commonly appearing in the solutions of this type of model. To illustrate the performance of the method, several different test cases and benchmark problems have been run, showing the high-order of the scheme and the good agreement when compared to reference solutions computed with other well-known methods.

\end{abstract}

\section{Introduction}
\label{sec:introduction} 

Multi-phase flows are ubiquitous in nature and engineering applications. The simplest flow of a liquid with free-surface, such as a flowing river or a falling raindrop, already involves both the dynamics of the liquid phase and the surrounding air and can consequently be considered a two-phase flow. The application range is obviously much larger and includes, for example, and without pretending to be exhaustive, bubbly liquids and sprays, water flow with sediment transport, mist flows, two-phase flows with phase change as used in modern 3D printers, compressible multi-phase flows in internal combustion engines, flows in the paper, steel and food industry, etc. Such applications have motivated extensive efforts to study and develop multi-phase flow models that describe such phenomena with respect to physical principles and thermodynamics. 

One particular aspect of multi-phase flows is that they involve moving interfaces between different bulk phases. From a physics perspective, the nature of the interface, that is, whether it is sharp or diffuse, has been subject to debate since the times of Rayleigh and Laplace. While both approaches have their advantages, it is admitted that diffuse interface approaches, in general, offer more flexibility, especially in the presence of strong deformations and topology changes in the interface geometry. In contrast, the sharp interface approach allows for a more rigorous treatment of the thermodynamics at the interface. In this paper, we are interested in diffuse interface approaches. Although a thorough review of diffuse interface approaches is beyond the scope of this work, the reader is, for example, referred  to~\cite{Anderson1998,Kapila2001,SaurelAbgrall1999MultiphaseGod,AbgrallSaurel2003,FavrieGavrilyukSaurel2009,FavrieGavrilyuk2012,NdanouFavrieGavrilyuk2015,Barton2019,DeLorenzo2018,Gaburro2018diffuse,Busto2020ADERGPR,Kemm2020,ReAbgrall2022,vanderWaals1979} and references therein. For the incompressible case, see also \cite{Casulli2014VOF,FerrariDumbser2021SIVOF,FerrariDumbser2023SIVOF2} and related work.  

Up to now no universally accepted mathematical model exists for the whole range of different compressible multi-phase flows. Among the large number of different models that can be found in the literature, the one proposed by Baer and Nunziato~\cite{BaerNunziato1986} is one of the most widespread, see, for example, \cite{AndrianovWarnecke2004,Bdzil1999,GavrilyukSaurel2002,SaurelAbgrall1999MultiphaseGod}, as well as the well-known Kapila model \cite{Kapila2001}, which can be obtained from the Baer--Nunziato model in the stiff relaxation limit. The governing equations of the Baer--Nunziato (BN) model form a hyperbolic system because the associated eigenvalues are all real, and there exists a set of linearly independent eigenvectors. However, some of its equations cannot be written in conservative form, which makes it difficult to deal with the appearance of discontinuities and the development of high-order numerical methods. Moreover, the BN model was modified by different authors in the literature, see, for example, Saurel and Abgrall~\cite{SaurelAbgrall1999MultiphaseGod}, whose modification describes multi-phase mixtures and interface problems between pure compressible materials. An alternative two-phase flow model that is fully conservative is the one forwarded by Scannapieco and Cheng, see \cite{ScannapiecoCheng2002}. 

In addition to the previously mentioned models, there exists another class of models for compressible multi-phase flows that originates from the theory of Godunov and Romenski on symmetric hyperbolic and thermodynamically compatible systems \cite{Godunov1961,Romenski1998,Romenski2001,Godunov2003} and which was first introduced in \cite{Romenski2004,Romenski2007TwoPhase,Romenski2010TwoPhase,RomenskiBelozerov2016}. In this paper, we are therefore interested in the discretization of the SHTC model of \textit{barotropic} compressible two-phase flows of Romenski \textit{et al.} with different phase velocities and phase pressures in multiple space dimensions. The model consists of a first-order symmetric hyperbolic and thermodynamically compatible (SHTC) system of equations~\cite{Romenski2007TwoPhase,Romenski2010TwoPhase,RomenskiBelozerov2016}. In the previous references, it was also shown that the SHTC equations, which are written in conservative form, can be converted to the form of a BN-type model, where additional differential terms appear in the momentum equations, which were not included in the original BN model and which describe the so-called lift forces.  Recently, an exact solution for the corresponding Riemann problem in the barotropic case was found in~\cite{Thein2022}, and an all-Mach number flow solver was developed in~\cite{Lukacova2022}. 

The model proposed by Romenski \textit{et al.} is strongly hyperbolic in one space dimension, and its homogeneous part without algebraic source terms is endowed with a curl involution on the relative velocity field. However, as we will show in this paper, the original model is only \text{weakly hyperbolic} in multiple space dimensions. In order to restore strong hyperbolicity, two different strategies can be followed: the first one consists in using the hyperbolic generalized Lagrangian multiplier (GLM) curl-cleaning methodology introduced in~\cite{Dumbser2020GLM,Chiocchetti2021SHTCSurfaceTension,Busto2021HyperbolicDispersion,Dhaouadi2022NSK}, which is a natural extension of the original ideas presented by Munz \textit{et al.} in~\cite{Munz2000,Dedner2002} on hyperbolic GLM divergence cleaning for the Maxwell and MHD equations, which contain the well-known divergence-free condition of the magnetic field. Curl constraints can also be found in many other first-order hyperbolic models, such as hyperbolic models for surface tension~\cite{Schmidmayer2017,Chiocchetti2021SHTCSurfaceTension,Chiocchetti2023}, first-order hyperbolic reformulations of the Navier--Stokes--Korteweg equations based on an augmented Lagrangian approach~\cite{Dhaouadi2022NSK,Dhaouadi2023expspnsk}, or first-order reductions of the Einstein field equations of general relativity \cite{Dumbser2020GLM}. The second methodology to recover the hyperbolicity involves the use of some extra terms in the momentum equation that symmetrize the system and is therefore directly based on the theory of Symmetric Hyperbolic and Thermodynamically Compatible (SHTC) systems, following~\cite{Powell1994MHD1,Powell1997MHD,Powell1999MHD,Godunov1972MHD}. In both cases, we will show that the system in the multidimensional case becomes again strongly hyperbolic. 

In the setting of this model, the solutions are often discontinuous in space, that is, solutions to Riemann problems. Therefore, it is necessary to consider an approach that is robust and accurate even in the presence of shock waves or discontinuities. To address the sharp gradients of the numerical solutions, a high-order ADER Discontinuous Galerkin (DG) finite element framework with \aposteriori sub-cell finite volume (FV) limiter is considered, see~\cite{Dumbser2008DG,Dumbser2014,Zanotti2015} for further details. The proposed method is high-order in space and time thanks to the use of the ADER approach of Toro and Titarev~\cite{Toro2002,Titarev2005,Toro2006}. To deal with spurious oscillations that may appear in the presence of discontinuities or shock waves, it makes use of an \aposteriori subcell finite volume limiter for high-order fully discrete one-step ADER-DG schemes presented in~\cite{Dumbser2014,Zanotti2015}, which follows the MOOD approach of Clain and Loub\`ere~\cite{Clain2011,Diot2012,Diot2013}. 

This paper is organized as follows. Section~\ref{sec:gov_eq} recalls the set of governing partial differential equations and shows the EOS that we will use in this work. In Section~\ref{sec:hyperbolicity}, the hyperbolicity of the system is studied in the multidimensional case, showing that the original model is only weakly hyperbolic, and two different strategies to recover strong hyperbolicity are presented. Section~\ref{sec:dg.method} introduces the high-order ADER discontinuous Galerkin scheme used in this paper to solve the model numerically. Section~\ref{sec:num_results} shows the results of several test cases and benchmark problems computed using the proposed method. Finally, Section~\ref{sec:conclusion} concludes the paper and summarizes the contributions of this work as well as its possible future extensions.

\section{Governing equations}\label{sec:gov_eq}

\noindent The governing equations of the barotropic compressible two-velocity, two-pressure two-fluid model of Romenski, written in terms of the specific total energy potential $E=E(\alpha_1,c_1,\rho,w^k)$ are given by
\begin{subequations}\label{eq:barotropic_ini}
	\begin{align}
		&\frac{\der \rho \alpha_1}{\dt}+ \frac{\der \rho\alpha_1 \vel^k}{\dx{k}}=-\frac{\rho}{\tau}E_{\alpha_1}, \label{eq:volume_frac_ini}\\
		&\frac{\der\rho c_1}{\dt}+\frac{\der(\rho c_1 \vel^k + \rho E_{\urel^k})}{\dx{k}}=0,\label{eq:mass_frac_ini}\\
		&\frac{\der\rho}{\dt}+\frac{\der\rho \vel^k}{\dx{k}}=0,\label{eq:mass_ini}\\
		&\frac{\der\rho \vel^i}{\dt}+ \frac{\der(\rho \vel^i \vel^k+p\delta_{ik}+\rho \urel^i E_{\urel^k} )}{\dx{k}}=g^i \rho,\label{eq:momentum_ini}\\
		&\frac{\der \urel^k}{\dt}+\frac{\der (\urel^l \vel^l + E_{c_1})}{\dx{k}}+ \vel^l\left(\frac{\der \urel^k}{\dx{l}}-\frac{\der \urel^l}{\dx{k}} \right)=-\zeta E_{\urel_k},
		\label{eq:relative_vel_ini}
	\end{align}
\end{subequations}                                                                                         
where $\alpha_1$ and $\alpha_2$ are the volume fractions of the first and second phases, verifying the saturation constraint $\alpha_1 + \alpha_2 = 1$, $\rho_1$ and $\rho_2$ are the mass densities of the first and second phase,  $\rho=\alpha_1\rho_1+\alpha_2\rho_2$ is the mixture mass density,  $c_j=\alpha_j\rho_j / \rho$ with $j=1,2$ are the mass fractions of phase $j$ which satisfy $c_1+c_2=1$, and $\bg=(g^1,g^2,g^3)^T$ is the gravity acceleration. If $\bvel_1=(\vel_1^1,\vel_1^2,\vel_1^3)^T$ and $\bvel_2=(\vel_2^1,\vel_2^2,\vel_2^3)^T$ are the velocity vectors of each phase, then the mixture velocity $\bvel=(\vel^1,\vel^2,\vel^3)^T$ is computed as $\vel^k=c_1\vel^k_1+c_2\vel^k_2$ and the relative phase velocity, which is a primary evolution quantity in this model, $  $is given by  $\burel=(\urel^1,\urel^2,\urel^3)^T$ as $\urel^k = \vel^k_1-\vel^k_2$. The source terms of equations \eqref{eq:volume_frac_ini} and~\eqref{eq:relative_vel_ini} are proportional to thermodynamic forces, with $\tau$ the rate of pressure relaxation and $\zeta$ the inter-phase friction coefficient. The PDE system~\eqref{eq:barotropic_ini} is formed by two conservation laws for the volume and mass fractions, Eq.~\eqref{eq:volume_frac_ini} and \eqref{eq:mass_frac_ini}, respectively, conservation of total mass, Eq.~\eqref{eq:mass_ini}, conservation of the mixture momentum, Eq.~\eqref{eq:momentum_ini}, and a balance law for the relative velocity, Eq.~\eqref{eq:relative_vel_ini}. The algebraic source terms appearing in Eq.~\eqref{eq:volume_frac_ini} and~\eqref{eq:relative_vel_ini}, describe the interaction between the phases and are pressure relaxation and interfacial friction.
The mixture equation of state is given by, see~\cite{Thein2022}, 
\begin{equation*}
	E(\alpha_1,c_1,\rho,\urel^1, \urel^2, \urel^3) = e(\alpha1,c_1,\rho)+c_1c_2\frac{w_iw^i}{2},
\end{equation*}
where the specific internal energy of the mixture reads as
\begin{equation*}
	e(\alpha_1,c_1,\rho)=c_1e_1(\rho_1)+c_2e_2(\rho_2) = c_1e_1\left(\frac{c_1\rho}{\alpha_1}\right)+c_2e_2\left(\frac{c_2\rho}{\alpha_2} \right),
\end{equation*}
where $e_j(\rho_j)$ is the specific internal energy of the phase $j$. Since in this work an isentropic process is considered, the derivatives of $e$ can be computed as 
\begin{equation*}
	\frac{\partial e}{\partial \alpha_1} = \frac{p_2 -p_1}{\rho}, \qquad
	\frac{\partial e}{\partial c_1} = h_1(\rho_1) - h_2(\rho_2), \qquad
	\frac{\partial e}{\partial \rho} = \frac{p}{\rho^2},
\end{equation*}
where $h_j(\rho_j) = e_j(\rho_j)+\frac{p_j(\rho_j)}{\rho_j}, \, j=1,2$ is the specific enthalpy of the phase $j$. Then, the derivatives of the specific total energy $E$ are given by
\begin{align}
	\frac{\partial E}{\partial \alpha_1} & = \frac{\partial e}{\partial \alpha_1} = \frac{p_2 -p_1}{\rho}, \qquad \frac{\partial E}{\partial \rho} = \frac{\partial e}{\partial \rho}= \frac{p}{\rho^2}, \qquad 	\frac{\partial E}{\partial \urel^i} = c_1 (1-c_1) \urel^i,\label{eq:der_specific_energy1}\\
	\frac{\partial E}{\partial c_1} & = \frac{\partial e}{\partial c_1} + (1-2c_1)\frac{\urel_i\urel^i}{2} = h_1(\rho_1) - h_2(\rho_2) + (1-2c_1)\frac{\urel_i\urel^i}{2}.
	\label{eq:der_specific_energy2}
\end{align}

Considering the computations in~\cite{Thein2022}, and taking into account~\eqref{eq:der_specific_energy1} and~\eqref{eq:der_specific_energy2}, the following identities are obtained
\begin{align*}
	&\rho c_1 u^k + \rho E_{w^k}= \rho c_1 u_1^k, \qquad 
	 \rho c_2 u^k - \rho E_{w^k}= \rho c_2 u_2^k,\\
	&\rho u^i u^k+p\delta_{ik}+\rho w^i E_{w^k}= \alpha_1 \rho_1 u_1^i u_1^k + \alpha_2 \rho_2 u_2^i u_2^k + p \delta_{ik},\\
	&(u_1^l-u_2^l) u^l + E_{c_1}= \frac{1}{2}(u_1^l)^2 -\frac{1}{2} (u_2^l)^2 + h_1-h_2.
\end{align*}
The PDE system~\eqref{eq:barotropic_ini} can now be rewritten in the following form, which is more convenient for numerical discretization:  
\begin{subequations}\label{eq:cons_barotropic_system}
	\begin{align}
		&\frac{\dalpha{1}}{\dt}+u^k \frac{\dalpha{1}}{\dx{k}}=-\frac{p_2-p_1}{\tau \rho},\label{eq:cons_volume_frac}\\
		&\frac{\der\alpha_1\rho_1}{\dt}+\frac{\der(\alpha_1\rho_1 u_1^k)}{\dx{k}}=0,\label{eq:cons_mass_density1}\\
		&\frac{\der\alpha_2\rho_2}{\dt}+\frac{\der(\alpha_2\rho_2 u_2^k)}{\dx{k}}=0,\label{eq:cons_mass_density2}\\
		&\frac{\der\rho u^i}{\dt}+ \frac{\der(\alpha_1\rho_1 u_1^i u_1^k+ \alpha_2\rho_2 u_2^i u_2^k + p\delta_{ik} )}{\dx{k}}=g_i \rho,\label{eq:cons_momentum}\\
		&\frac{\der w^k}{\dt}+\frac{\der}{\dx{k}} \left( \frac{1}{2}(u_1^l)^2-\frac{1}{2}(u_2^l)^2 + h_1-h_2\right)+ u^l\left(\frac{\der w^k}{\dx{l}}-\frac{\der w^l}{\dx{k}} \right)=-\zeta c_1c_2w^k.\label{eq:cons_relative_velocity}
	\end{align}
\end{subequations}
The derivation of the PDE system~\eqref{eq:cons_barotropic_system} is based on the principles of thermodynamically compatible systems \cite{Godunov2003}, and it consists of nine equations: the balance law for the volume fraction, the conservation laws of the two-phase masses, conservation of mixture momentum and the balance law for the relative velocity field. Note that in the absence of source terms in \eqref{eq:cons_relative_velocity}, i.e., for $\zeta = 0$, the relative velocity is curl-free in the sense
\begin{equation*}
	\frac{\der w^k}{\dx{l}}-\frac{\der w^l}{\dx{k}}  = 0.   
\end{equation*} 

\subsection{Equation of state (EOS)}\label{sec:eos}
In order to close the two-phase model \eqref{eq:cons_volume_frac}-\eqref{eq:cons_relative_velocity}, it is necessary to define an equation of state for each phase: throughout this paper, we will make use either of the ideal gas law or of the stiffened gas equation of state, which will be used to model a liquid phase.
For ideal gases, the EOS is defined as 
\begin{equation}
E(\rho)=\frac{c_0^2}{\gamma (\gamma-1)} \qquad \mbox{with } c_0^2=\gamma \rho^{\gamma-1}e^{s/c_v},
\label{eq:energy_ideal}
\end{equation}
where $\gamma$ is the adiabatic index or the ratio of specific heats, $c_0$ is the adiabatic sound speed, $s$ is the specific entropy (which in our case will be constant), $c_v$ is the specific heat capacity at constant volume, and the pressure is given by
\begin{equation}
p(\rho) = \rho^2 E_{\rho} = \rho ^\gamma e^{s/c_v} = (\gamma-1)\rho E.
\label{eq:pressure_ideal}
\end{equation}
For stiffened gases, the EOS reads as
\begin{equation}
E(\rho)=\frac{c_0^2}{\gamma (\gamma-1)}\left(\frac{\rho}{\rho_0}\right)^{\gamma-1}e^{s/c_v}+\frac{\rho_0 c_0^2 - \gamma p_0}{\gamma \rho},
\label{eq:energy_stiffened}
\end{equation}
where $\rho_0$ and $p_0$ are the reference density and pressure, respectively, and $c_0$ is a constant reference sound speed. In this case, the pressure is computed as
\begin{equation}
p(\rho) = \rho^2 E_{\rho} = \frac{c_0^2 \rho_0}{\gamma} \left(\frac{\rho}{\rho_0} \right)^\gamma e^{s/c_v}- \frac{c_0^2 \rho_0 - \gamma p_0}{\gamma}.
\label{eq:pressure_stiffened}
\end{equation}

\subsection{Hyperbolicity analysis}\label{sec:hyperbolicity}
In this section, we will study the hyperbolicity of the model~\eqref{eq:cons_volume_frac}-\eqref{eq:cons_relative_velocity}. The one-dimensional case was already addressed in~\cite{Thein2022}, and here we briefly recall some results for the one-dimensional case before moving to the more general multidimensional case. In particular, we will prove that the original system \eqref{eq:cons_volume_frac}-\eqref{eq:cons_relative_velocity} is only \textit{weakly hyperbolic} in the multidimensional case and show how the strong hyperbolicity can be restored considering two different strategies: on the one hand, using an extension of the hyperbolic Generalized Lagrangian Multiplier (GLM) curl-cleaning approach and on the other hand, modifying the original system of governing equations by adding the symmetrizing terms that allow rewriting the model in symmetric hyperbolic form, which is natural within the framework of SHTC equations. Defining the vectors of conserved and primitive variables as $\bQ=(\alpha_1,\alpha_1\rho_1,\alpha_2\rho_2, \rho \bvel^T, \burel^T)^T$ and $\bV =(\alpha_1,\rho_1,\rho_2,\bvel_1^T,\bvel_2^T)^T$, respectively, the PDE system~\eqref{eq:cons_barotropic_system} can be written as
\begin{equation*}
	\partial_t \bQ + \nabla \cdot \bF(\bQ) + \bB{}(\bQ)\cdot\nabla\bQ =\bS(\bQ),
\end{equation*}
where $\bS(\bQ)$ contains the source terms, $\bF(\bQ)$ is the nonlinear flux tensor and $\bB{}(\bQ)\cdot\nabla\bQ$ contains the non-conservative terms. Then, the quasi-linear form of the PDE in terms of the conserved variables $\bQ$ is  given by
\begin{equation*}
	\partial_t \bQ + \bA(\bQ)\cdot\nabla\bQ = \bS(\bQ),
\end{equation*}
where $\bA(\bQ) = \ds{\frac{\partial \bF}{\partial \bQ} + \bB{}}$. If the vector of primitive variables $\bV$ is considered, the system can be written as
\begin{equation}
	\partial_t \bV + \bC(\bV)\cdot\nabla\bV = \bcalS(\bV),
	\label{eq:PDEcompactV}
\end{equation}
where $\bC(\bV) = \ds{\frac{\partial \bV}{\partial \bQ} \frac{\partial \bF}{\partial \bV} + \frac{\partial \bV}{\partial \bQ} \bB{}\frac{\partial \bQ}{\partial \bV}}$ and $\bcalS(\bV)=\ds{\frac{\partial \bV}{\partial \bQ} \, \bS(\bQ)}$. 
Throughout this section, for the sake of readability and since it does not contribute anything to the study of the hyperbolicity of the system, the source terms are set to 0.

\subsubsection{One-dimensional case}
The hyperbolicity analysis of the system~\eqref{eq:cons_barotropic_system} in 1D has been done in detail in~\cite{Thein2022}, so here, only a summary of the main points will be made that will be useful for a better understanding of the multidimensional case. If $\vel = c_1 \vel_1+ c_2 \vel_2$ is the mixture velocity, with $\vel_j,\, j=1,2$ the velocity of the phase $j$ and $\urel =\vel_1-\vel_2$ is the relative velocity, the one-dimensional system results
\begin{subequations}\label{eq:cons_barotropic_system_1D}
	\begin{align}
		&\frac{\dalpha{1}}{\dt}+\vel \frac{\dalpha{1}}{\dx{}}=0,\label{eq:cons_volume_frac_1D}\\
		&\frac{\der\alpha_1\rho_1}{\dt}+\frac{\der(\alpha_1\rho_1 \vel_1)}{\dx{}}=0,\label{eq:cons_mass_density1_1D}\\
		&\frac{\der\alpha_2\rho_2}{\dt}+\frac{\der(\alpha_2\rho_2 \vel_2)}{\dx{}}=0,\label{eq:cons_mass_density2_1D}\\
		&\frac{\der\rho \vel}{\dt}+ \frac{\der(\alpha_1\rho_1 (\vel_1)^2+ \alpha_2\rho_2 (\vel_2) + p )}{\dx{}}=0,\label{eq:cons_momentum_1D}\\
		&\frac{\der \urel}{\dt}+\frac{\der}{\dx{}} \left( \frac{1}{2}(\vel_1)^2-\frac{1}{2}(\vel_2)^2 + h_1-h_2\right)=0.\label{eq:cons_relative_velocity_1D}
	\end{align}
\end{subequations}
Then, the matrix $\bC$ in~\eqref{eq:PDEcompactV} is given by
\begin{equation*}	
	\bC(\bV) = \left(
	\arraycolsep=4.pt\def\arraystretch{1.4}
	\begin{array}{ccccc}
		\vel & 0 & 0 & 0 & 0 \\
		\frac{ \rho_1}{\alpha_1} \left(\vel_1-\vel\right) & u_1 & 0 & \rho_1 & 0 \\
		\frac{ \rho_2}{\alpha_2} \left(\vel-\vel_2\right) & 0 & u_2 & 0 & \rho_2 \\
		\beta & \frac{a_1^2}{\rho_1} & 0 & \vel_1 & 0 \\
		\beta & 0 & \frac{a_2^2}{\rho_2} & 0 & \vel_2 \\
	\end{array}
	\right),
	\label{eq:Jacobian_1D}
\end{equation*}
with $a_j$ the sound speed of phase $j$, that is defined as $a^2_j=\rho_j\left(\frac{\partial h_j}{\partial \rho_j}\right), \, j=1,2$ and $\beta=\frac{1}{\rho}\left(\frac{a_1^2\rho_1}{\gamma_1}- \frac{a_2^2\rho_2}{\gamma_2}\right) = \frac{1}{\rho}\left(p_1- p_2\right)$. It is easy to show that $\bC(\bV)$ admits five eigenvalues, whose expressions are given hereafter 
\begin{equation*}
	\lambda_1 = u_1+a_1, \,\,
	\lambda_2 = u_1-a_1, \,\,
	\lambda_3 = u, \,\, 
	\lambda_4 = u_2+a_2,\,\,
	\lambda_5 = u_2-a_2. 
\end{equation*}
In this case, all five eigenvalues are real, and a complete set of five linearly independent eigenvectors exists, which means that the system in one space dimension is \textit{strongly hyperbolic}, see~\cite{Thein2022} for further details about the eigenvalues and eigenvectors in the one-dimensional case.

\subsubsection{Multidimensional case}
Now we are in the position to study the multidimensional case. In the following, we use of the property of \textit{rotational invariance} of Newtonian mechanics, hence it is enough to consider the matrix $\bC(\bV)$ only in the $x-$direction and not all possible space directions. The system under consideration is~\eqref{eq:cons_barotropic_system} and, in this case, the matrix $\bC(\bV)$ in the $x-$direction is given by
\begin{equation*}
	\bC(\bV) = \left(
	\arraycolsep=3.pt\def\arraystretch{1.2}
	\begin{array}{ccccccccc}
		u^1 & 0 & 0 & 0 & 0 & 0 & 0 & 0 & 0  \\
		\frac{ \rho_1}{\alpha_1}\vartheta_2^1 & u^1_1 & 0 & \rho_1 & 0 & 0 & 0 & 0 & 0  \\
		\frac{ \rho_2}{\alpha_2}\vartheta_1^1& 0 & u_2^1 & 0 & 0 & 0 & \rho_2 & 0 & 0 \\
		\beta & \frac{a_1^2}{\rho_1} & 0 & u_1^1 & \vartheta_{22}^2 & \vartheta_{22}^3& 0 & \vartheta_{12}^2& \vartheta_{12}^3\\
		0 & 0 & 0 & 0 & \xi_2 & 0 & 0 & -\vartheta_{12}^1 & 0  \\
		0 & 0 & 0 & 0 & 0 & \xi_2 & 0 & 0 & -\vartheta_{12}^1\\
		\beta & 0 & \frac{a_2^2}{\rho_2} & 0 & -\vartheta_{12}^2 & -\vartheta_{12}^3 & u_2^1 & -\vartheta_{11}^2 & -\vartheta_{11}^3 \\
		0 & 0 & 0 & 0 & \vartheta_{12}^1 & 0 & 0 & \xi_1 & 0  \\
		0 & 0 & 0 & 0 & 0 & \vartheta_{12}^1 & 0 & 0 & \xi_1  \\
	\end{array}
	\right).
\end{equation*}
As in the previous section, we make use of the following auxiliary variables to ease the expressions in the matrix
\begin{equation}
\left.
\begin{array}{l}
\beta= \frac{1}{\rho}(p_1-p_2),\\[0.2cm]
\vartheta_j^i=c_j w^i,\,\, \vartheta_{jk}^i=c_j c_k w^i,\,\, i=1,2,3,\,\, j=1,2,\,\, k=1,2,\\[0.2cm]
\xi_1 = c_1u^1+c_2u_2^1 = u_2^1 + c_1^2 \urel^1 = u_2^1 + c_1 \vartheta_1^1, \\[0.2cm]
\xi_2 = c_1u_1^1+c_2u^1 = u_1^1 - c_2^2 \urel^1= u_1^1 - c_2 \vartheta_2^1.
\end{array}\right.
\label{eq:auxiliary_var_jacobian}
\end{equation}
The matrix $\bC$ admits 9 eigenvalues $\lambda_{1-9}$ that are given by
\begin{equation*}
\lambda_1 = u_1^1-a_1, \,\,
\lambda_2 = u_1^1+a_1, \,\, 
\lambda_{3} = u_2^1-a_2,\,\,
\lambda_{4} = u_2^1+a_2, \,\, 
\lambda_{5-9} = u^1,   
\end{equation*}
where the sound speed of each phase $a_j$ again is defined as $a^2_j=\rho_i\left(\frac{\partial h_j}{\partial \rho_j}\right),\,\,j=1,2$. The eigenvalues are all real. To prove whether the system is weakly or strongly hyperbolic, it is necessary to compute the associated eigenvectors. The right eigenvectors are the columns in the matrix below and are given in the same order as the eigenvalues:  
\begin{equation*}
	\bR_{\scriptscriptstyle{1-4}}=\left(
	\arraycolsep=1.8pt\def\arraystretch{1.5}
	\begin{array}{cccc}
		0 & 0 & 0 & 0\\
	    {-\frac{\rho_1}{a_1}} & {\frac{\rho_1}{a_1}} & 0 & 0\\
		0 & 0 &	{-\frac{\rho_2}{a_2}} & {\frac{\rho_2}{a_2}}\\
		1 & 1 & 0 & 0\\
		0 & 0 & 0 & 0\\
		0 & 0 & 0 & 0\\
		0 & 0 & 1 & 1\\
		0 & 0 & 0 & 0\\
		0 & 0 & 0 & 0\\
	\end{array}
	\right), \,\,
	\bR_{\scriptscriptstyle{5-7}}=\left(
	\arraycolsep=1.8pt\def\arraystretch{1.5} 
	\begin{array}{ccc}
		{-\frac{\alpha_1\,\alpha_2\,Z_{2}}{\eta_2\,\vartheta_{1}^1}} & 0 & 0 \\
		{\frac{\rho_1\,\zeta_1\,Z_{2}}{\eta_2\,\vartheta_{1}^1 Z_{1}}} & {-\frac{ \rho_1 \,\vartheta_{2}^3}{Z_{1}}} & {-\frac{ \rho_1\,\vartheta_{2}^2}{Z_{1}}} \\
		{\frac{\rho_2\,\zeta_2}{\eta_2\,\vartheta_{1}^1}} & {\frac{\rho_2\, \vartheta_{1}^3}{Z_{2}}} & {\frac{\rho_2\,\vartheta_{1}^2}{Z_{2}}} \\
		{-\frac{c_2\,\eta_1\,\vartheta_{2}^1}{c_1\,\eta_2\,Z_{1}}}  & {\frac{\vartheta_{2}^1\,\vartheta_{2}^3}{Z_{1}}} & {\frac{\vartheta_{2}^1\, \vartheta_{2}^2}{Z_{1}}} \\
		0 & 0 & 1 \\
		0 & 1 & 0 \\
		1 & {\frac{\vartheta_{1}^1\,\vartheta_{1}^3}{Z_{2}}} & {\frac{\vartheta_{1}^1\, \vartheta_{1}^2}{Z_{2}}} \\
		0 & 0 & 1 \\
		0 & 1 & 0 \\
	\end{array}
	\right),
\end{equation*}
where the auxiliary variables that are used to ease the notation are defined as
\begin{align*}
		Z_{1}= a_1^2-(\vartheta_{2}^1)^2, \quad 
		Z_{2}=a_2^2-(\vartheta_{1}^1)^2,  \quad  
		\eta_1  = \alpha_2(a_1^2-\alpha_1\beta),  \nonumber \\ 
		\eta_2  = \alpha_1(a_2^2+\alpha_2\beta), \quad 
		\zeta_1 = \alpha_2\left(\alpha_1\beta-(\vartheta_2^1)^2\right), \quad 
		\zeta_2 = \alpha_1\left(\alpha_2\beta+(\vartheta_1^1)^2\right).
\end{align*}  
All eigenvalues are real, but two eigenvectors are missing (in fact, if we are working in dimension $d$ there are $d-1$ missing eigenvectors), so the system is only \textit{weakly hyperbolic}. Hereafter we will show two different methodologies that can be used to restore the strong hyperbolicity of the model.

\paragraph{Generalized Lagrangian multiplier (GLM) curl-cleaning approach}
To restore the strong hyperbolicity and following the same strategy that can be found in~\cite{Dumbser2020GLM,Chiocchetti2021SHTCSurfaceTension,Busto2021HyperbolicDispersion,Dhaouadi2022NSK}, we make use of the GLM curl-cleaning technique, where an evolution equation for a curl-cleaning field $\psi^k$ is added to the system~\eqref{eq:cons_barotropic_system}. Using the abbreviation $ \delta_{12} = \frac{1}{2}(u_1^l)^2-\frac{1}{2}(u_2^l)^2 + h_1-h_2$ and denoting the curl-cleaning speed by $a_{\bpsi}$ this equation is coupled with~\eqref{eq:cons_relative_velocity} via a Maxwell-type sub-system as follows:

\begin{eqnarray*}
	\ds{\frac{\der w^k}{\dt}+\frac{\der  \delta_{12}}{\dx{k}} + u^l\left(\frac{\der w^k}{\dx{l}}-\frac{\der w^l}{\dx{k}} \right)} \ds{+ a_{\bpsi}  \varepsilon_{klm} \frac{\der \psi^m}{\dx{l}}=0,} \\ 
	\ds{\frac{\der\psi^k}{\dt} + u^j\frac{\der \psi^k}{\dx{j}} -a_{\bpsi} \varepsilon_{klm}\frac{\der \urel^m}{\dx{l}} = 0, }
\end{eqnarray*}
hence the augmented system with GLM curl-cleaning reads 
\begin{subequations}\label{eq:GLM_barotropic_system}
	\begin{align}
		&\frac{\dalpha{1}}{\dt}+u^k \frac{\dalpha{1}}{\dx{k}}=0,\label{eq:GLM_volume_frac}\\
		&\frac{\der\alpha_1\rho_1}{\dt}+\frac{\der(\alpha_1\rho_1 u_1^k)}{\dx{k}}=0,\label{eq:GLM_mass_density1}\\
		&\frac{\der\alpha_2\rho_2}{\dt}+\frac{\der(\alpha_2\rho_2 u_2^k)}{\dx{k}}=0,\label{eq:GLM_mass_density2}\\
		&\frac{\der\rho u^i}{\dt}+ \frac{\der(\alpha_1\rho_1 u_1^i u_1^k+ \alpha_2\rho_2 u_2^i u_2^k + p\delta_{ik} )}{\dx{k}}=0 ,\label{eq:GLM_momentum}\\
		&\frac{\der w^k}{\dt}+\frac{\der \delta_{12}}{\dx{k}} + u^l\left(\frac{\der w^k}{\dx{l}}-\frac{\der w^l}{\dx{k}} \right) + a_{\bpsi}  \varepsilon_{klm} \frac{\der \psi^m}{\dx{l}}=0,\label{eq:GLM_relative_velocity}\\
		&\frac{\der\psi^k}{\dt} + u^j\frac{\der \psi^k}{\dx{j}} -a_{\bpsi} \varepsilon_{klm}\frac{\der \urel^m}{\dx{l}} = 0, \label{eq:GLM_psi}
	\end{align}
\end{subequations}
where $\bpsi=(\psi^1,\psi^2,\psi^3)$ is the cleaning field, $\boldsymbol{\varepsilon}=\varepsilon_{klm}$ is the Levi–Civita tensor and $a_{\bpsi}$ is the curl-cleaning speed. Once the curl-cleaning field has been added, we can compute the eigenvalues and eigenvectors to check the hyperbolicity of the augmented GLM system. The matrix $\bC$ in $x-$direction can be written as
\begin{equation*}
\bC(\bV)=\left(
	\arraycolsep=2.5pt\def\arraystretch{1.2}
	\begin{array}{cccccccccccc}
		u^1 & 0 & 0 & 0 & 0 & 0 & 0 & 0 & 0 & 0 & 0 & 0 \\
		\frac{ \rho_1}{\alpha_1}\vartheta_2^1 & u^1_1 & 0 & \rho_1 & 0 & 0 & 0 & 0 & 0 & 0 & 0 & 0 \\
		\frac{ \rho_2}{\alpha_2}\vartheta_1^1 & 0 & u_2^1 & 0 & 0 & 0 & \rho_2 & 0 & 0 & 0 & 0 & 0 \\
		\beta & \frac{a_1^2}{\rho_1} & 0 & u_1^1 & \vartheta_{22}^2 & \vartheta_{22}^3& 0 & \vartheta_{12}^2& \vartheta_{12}^3 & 0 & 0 & 0 \\
		0 & 0 & 0 & 0 & \xi_2 & 0 & 0 & -\vartheta_{12}^1 & 0 & 0 & 0 & -c_2 a_{\bpsi} \\
		0 & 0 & 0 & 0 & 0 & \xi_2 & 0 & 0 & -\vartheta_{12}^1 & 0 & c_2 a_{\bpsi} & 0 \\
		\beta & 0 & \frac{a_2^2}{\rho_2}  & 0 & -\vartheta_{12}^2 & -\vartheta_{12}^3 & u_2^1 & -\vartheta_{11}^2 & -\vartheta_{11}^3 & 0 & 0 & 0 \\
		0 & 0 & 0 & 0 & \vartheta_{12}^1 & 0 & 0 & \xi_1 & 0 & 0 & 0 & c_1 a_{\bpsi}  \\
		0 & 0 & 0 & 0 & 0 & \vartheta_{12}^1 & 0 & 0 & \xi_1 & 0 & -c_1 a_{\bpsi}  & 0 \\
		0 & 0 & 0 & 0 & 0 & 0 & 0 & 0 & 0 & u^1 & 0 & 0 \\
		0 & 0 & 0 & 0 & 0 & a_{\bpsi} & 0 & 0 & -a_{\bpsi} & 0 & u^1 & 0 \\
		0 & 0 & 0 & 0 & -a_{\bpsi} & 0 & 0 & a_{\bpsi} & 0 & 0 & 0 & u^1 \\
	\end{array}
	\right),
\end{equation*}
where, as in the previous case, we make use of the auxiliary variables~\eqref{eq:auxiliary_var_jacobian} to ease the notation of the matrix.
In this case, it is easy to compute the twelve eigenvalues of the matrix $\bC$ that are given by
\begin{align*}
	&\lambda_1 = u^1_1-a_1,    \quad  
	\lambda_2 = u^1_1+a_1,    \quad 
	\lambda_{3} = u^1_2-a_2, \quad  
	\lambda_{4} = u^1_2+a_2, \\ 
	&\lambda_{5-8} = u^1,    \quad 
	\lambda_{9-10} = u^1-a_{\bpsi},  \quad  
	\lambda_{11-12} = u^1+a_{\bpsi},
\end{align*}
with $a_j$ the sound speed of phase $j$, that is defined as $a^2_j=\rho_j\left(\frac{\partial h_j}{\partial \rho_j}\right), \, j=1,2$. Since the eigenvalues are real, to check if the system is weakly or strongly hyperbolic, it is necessary to compute the associated eigenvectors. Below, we will write the matrix that contains the right eigenvectors in columns. They are listed in the same order as the eigenvalues.
\begin{equation*}
	\bR_{\scriptscriptstyle{1-8}}=\left(
	\arraycolsep=2.2pt\def\arraystretch{1.5}
	\begin{array}{cccccccc}	 
		0 & 0 & 0 & 0 & 0 & {\frac{\alpha_1 \, \alpha_2 \,\vartheta_1^3}{\eta_2}} & {\frac{\alpha_1\, \alpha_2 \,\vartheta_1^2}{\eta_2}} & {-\frac{\alpha_1\, \alpha_2\, Z_{2}}{\eta_2\, \vartheta_1^1}} \\[0.2cm]
		{-\frac{\rho_1}{a_1}} & {\frac{\rho_1}{a_1}} & 0 & 0 & 0 & {-\frac{\varrho\, \vartheta_1^3}{\eta_2\,Z_{1}}} & {-\frac{\varrho \,\vartheta_1^2}{\eta_2\,Z_{1}}} & {\frac{\zeta_1\,\rho_1\, Z_{2}}{\eta_2\,\vartheta_1^1\,Z_{1}}} \\[0.2cm]
		0 & 0 & {-\frac{\rho_2}{a_2}} & {\frac{\rho_2}{a_2}} & 0 & {\frac{\alpha_1\, \rho_2\, \vartheta_1^3}{\eta_2}} & {\frac{\alpha_1\, \rho_2\, \vartheta_1^2}{\eta_2}} & {\frac{\zeta_2\, \rho_2}{\eta_2\,\vartheta_1^1}}\\[0.2cm]
		1 & 1 & 0 & 0 & 0 &{-\frac{\varepsilon\, \vartheta_2^1\, \vartheta_1^3}{\eta_2\,\,Z_{1}}} & {-\frac{ \varepsilon\, \vartheta_1^2 \vartheta_2^1}{\eta_2\,Z_{1}}} & {-\frac{c_2\,\eta_1 Z_{2}}{c_1\,\eta_2 Z_{1}}} \\	 
		0 & 0 & 0 & 0 & 0 & 0 & 1 & 0 \\
		0 & 0 & 0 & 0 & 0 & 1 & 0 & 0 \\
		0 & 0 & 1 & 1 & 0 & 0 & 0 & 1 \\
		0 & 0 & 0 & 0 & 0 & 0 & 1 & 0 \\
		0 & 0 & 0 & 0 & 0 & 1 & 0 & 0 \\
		0 & 0 & 0 & 0 & 1 & 0 & 0 & 0 \\
		0 & 0 & 0 & 0 & 0 & 0 & 0 & 0 \\
		0 & 0 & 0 & 0 & 0 & 0 & 0 & 0 \\
	\end{array}
	\right),
\end{equation*}

\begin{equation*}
\bR_{\scriptscriptstyle{9-12}}=\left(
\arraycolsep=1.8pt\def\arraystretch{1.5}
\begin{array}{cccc}
	0 & 0 & 0 & 0 \\
	-\frac{\rho_1 \tau_2^{+} \vartheta_2^2}{a_{\bpsi}\chi_2^{+}} & \frac{\rho_1 \tau_2^{+} \vartheta_2^3}{a_{\bpsi}\chi_2^{+}}&
	{\frac{\rho_1 \tau_2^{-} \vartheta_2^2}{a_{\bpsi}\chi_2^{-}}} & -\frac{\rho_1 \tau_2^{-} \vartheta_2^3}{a_{\bpsi}\chi_2^{-}}  \\[0.2cm]
	\frac{\rho_2 \tau_2^{+} \vartheta_1^2}{a_{\bpsi}\chi_1^{-}} & -\frac{\rho_2 \tau_2^{+} \vartheta_1^3}{a_{\bpsi}\chi_1^{-}}&
	-\frac{\rho_2 \tau_2^{-} \vartheta_1^2}{a_{\bpsi}\chi_1^{+}} & \frac{\rho_2 \tau_2^{-}	\vartheta_1^3}{a_{\bpsi}\chi_1^{+}}  \\[0.2cm]
	\frac{\delta_2^{+} \tau_2^{+} \vartheta_2^2}{a_{\bpsi}\chi_2^{+}} & -\frac{\delta_2^{+}	\tau_2^{+} \vartheta_2^3}{a_{\bpsi}\chi_2^{+}} &
	\frac{\delta_2^{-} \tau_2^{-} \vartheta_2^2}{a_{\bpsi}\chi_2^{-}} & -\frac{\delta_2^{-} \tau_2^{-} \vartheta_2^3}{a_{\bpsi}\chi_2^{-}} \\
	\frac{\upsilon_1^{-}}{a_{\bpsi}} & 0 & -\frac{\upsilon_1^{+}}{a_{\bpsi}} & 0 \\ 
	0 & -\frac{\upsilon_1^{-}}{a_{\bpsi}} & 0 & \frac{\upsilon_1^{+}}{a_{\bpsi}} \\[0.2cm] 
	\frac{\delta_1^{-} \tau_1^{-} \vartheta_1^2}{a_{\bpsi} \chi_1^{-}} & -\frac{\delta_1^{-} \tau_1^{-} \vartheta_1^3}{a_{\bpsi} \chi_1^{-}} &
	\frac{\delta_1^{+} \tau_1^{+} \vartheta_1^2}{a_{\bpsi} \chi_1^{+}} & -\frac{\delta_1^{+} \tau_1^{+}	\vartheta_1^3}{a_{\bpsi} \chi_1^{+}} \\
	-\frac{\upsilon_2^{+}}{a_{\bpsi}} & 0 & \frac{\upsilon_2^{-}}{a_{\bpsi}} & 0 \\
	0 & \frac{\upsilon_2^{+}}{a_{\bpsi}} & 0 & -\frac{\upsilon_2^{-}}{a_{\bpsi}} \\
	0 & 0 & 0 & 0 \\
	0 & 1 & 0 & 1 \\
	1 & 0 & 1 & 0 \\
\end{array}
\right),
\end{equation*}
where the auxiliary variables that have been used to write the twelve right eigenvectors more compactly are defined as
\begin{align*}
	& Z_{1}= a_1^2-(\vartheta_{2}^1)^2, \!\! \quad 
	  Z_{2}=a_2^2-(\vartheta_{1}^1)^2,  \!\! \quad 
  	\eta_1 = \alpha_2(a_1^2-\alpha_1\beta), \!\! \quad 
  	\eta_2 = \alpha_1(a_2^2+\alpha_2\beta), \\  
  	&\beta=\frac{p_1-p_2}{\rho}, \!\! \quad 
	\vartheta_j^i=c_j w^i, \!\! \quad
	\vartheta_{jk}^i=c_j c_k w^i, \!\! \quad 
	\xi_1 = u_2^1 + c_1 \vartheta_1^1, \!\! \quad 
	\xi_2 = u_1^1 - c_2 \vartheta_2^1, \\
	& \delta_j^{\pm}=a_{\bpsi}\pm \vartheta_j^1, \quad 
	\chi_1^{\pm}=(\delta_1^{\pm})^2-a_2^2, \quad  \chi_2^{\pm}=(\delta_2^{\pm})^2-a_1^2, \\  
	& \upsilon_1^{\pm}=c_2 \delta_1^{\pm}, \quad 
	  \upsilon_2^{\pm}=c_1 \delta_2^{\pm}, \quad  
	  \tau_1^{\pm}=\upsilon_1^{\pm}-c_1 a_{\bpsi}, \quad   \tau_2^{\pm}=\upsilon_2^{\pm}-c_2 a_{\bpsi}, \\  	  
	& \varrho = \alpha_2 \left(p_1+p_2(\gamma_2-1)-\rho_1 (\vartheta_2^1)^2\right), \quad 
	\varepsilon =\frac{\alpha_2}{\rho_1}\left(p_1(\gamma_1-1)-p_2(\gamma_2-1) \right).
\end{align*} 
Using the GLM curl-cleaning technique, we obtain twelve real eigenvalues and the corresponding twelve linearly independent right eigenvectors, hence the augmented system with GLM curl-cleaning~\eqref{eq:GLM_barotropic_system} is \textit{strongly hyperbolic}. 

\paragraph{Symmetrizing Godunov-Powell terms}
The second strategy used to recover the strong hyperbolicity of  system~\eqref{eq:cons_barotropic_system} is based on the theory of SHTC systems and consists in adding terms that are proportional to the curl involution to the momentum equation~\eqref{eq:cons_momentum} so that the system has real eigenvalues and a complete set of linearly independent eigenvectors, see \cite{Chiocchetti2021SHTCSurfaceTension}. The modified system reads
\begin{subequations}\label{eq:GodPowell_barotropic_ini}
	\begin{align}
		&\frac{\der \rho \alpha_1}{\dt}+ \frac{\der \rho\alpha_1 \vel^k}{\dx{k}}=0, \label{eq:GodPowell_volume_frac_ini}\\
		&\frac{\der\rho c_1}{\dt}+\frac{\der(\rho c_1 \vel^k + \rho E_{\urel^k})}{\dx{k}}=0,\label{eq:GodPowell_mass_frac_ini}\\
		&\frac{\der\rho}{\dt}+\frac{\der\rho \vel^k}{\dx{k}}=0,\label{eq:GodPowell_mass_ini}\\
		&\frac{\der\rho \vel^i}{\dt}+ \frac{\der(\rho \vel^i \vel^k+p\delta_{ik}+\rho \urel^i E_{\urel^k} )}{\dx{k}}+\rho E_{\urel^k}\left(\frac{\der w^k}{\dx{i}}-\frac{\der w^i}{\dx{k}}\right)=0,\label{eq:GodPowell_momentum_ini}\\
		&\frac{\der \urel^k}{\dt}+\frac{\der (\urel^l \vel^l + E_{c_1})}{\dx{k}}+ \vel^l\left(\frac{\der \urel^k}{\dx{l}}-\frac{\der \urel^l}{\dx{k}} \right)=0.
		\label{eq:GodPowell_relative_vel_ini}
	\end{align}
\end{subequations}
Using~\eqref{eq:der_specific_energy1}, system~\eqref{eq:GodPowell_barotropic_ini} results in 
\begin{subequations}\label{eq:GodPowell_barotropic_system}
	\begin{align}
		&\frac{\dalpha{1}}{\dt}+u^k \frac{\dalpha{1}}{\dx{k}}=0,\label{eq:GodPowell_volume_frac}\\
		&\frac{\der\alpha_1\rho_1}{\dt}+\frac{\der(\alpha_1\rho_1 u_1^k)}{\dx{k}}=0,\label{eq:GodPowell_mass_density1}\\
		&\frac{\der\alpha_2\rho_2}{\dt}+\frac{\der(\alpha_2\rho_2 u_2^k)}{\dx{k}}=0,\label{eq:GodPowell_mass_density2}\\
		&\frac{\der\rho u^i}{\dt}+ \frac{\der(\alpha_1\rho_1 u_1^i u_1^k+ \alpha_2\rho_2 u_2^i u_2^k + p\delta_{ik} )}{\dx{k}}+\rho c_1 c_2 \urel^k\left(\frac{\der w^k}{\dx{i}}-\frac{\der w^i}{\dx{k}}\right)=0 ,\label{eq:GodPowell_momentum}\\
		&\frac{\der w^k}{\dt}+\frac{\der}{\dx{k}} \left( \frac{1}{2}(u_1^l)^2-\frac{1}{2}(u_2^l)^2 + h_1-h_2\right)+ u^l\left(\frac{\der w^k}{\dx{l}}-\frac{\der w^l}{\dx{k}} \right)=0.\label{eq:GodPowell_relative_velocity}
	\end{align}
\end{subequations}
As in the previous case, the eigenvalues and eigenvectors are computed to analyze the hyperbolicity of the system~\eqref{eq:GodPowell_barotropic_system}. The matrix $\bC$ in the $x-$direction reads  
\begin{equation*}
	\bC(\bV) = \left(
	\arraycolsep=3.pt\def\arraystretch{1.2}
	\begin{array}{ccccccccc}
		u^1 & 0 & 0 & 0 & 0 & 0 & 0 & 0 & 0  \\
		\frac{ \rho_1}{\alpha_1}\vartheta_2^1 & u^1_1 & 0 & \rho_1 & 0 & 0 & 0 & 0 & 0  \\
		\frac{ \rho_2}{\alpha_2}\vartheta_1^1& 0 & u_2^1 & 0 & 0 & 0 & \rho_2 & 0 & 0 \\
		\beta & \frac{a_1^2}{\rho_1} & 0 & u_1^1 & \vartheta_{2}^2 & \vartheta_{2}^3& 0 & 0 & 0\\
		0 & 0 & 0 & 0 & u^1 & 0 & 0 & 0 & 0  \\
		0 & 0 & 0 & 0 & 0 & u^1 & 0 & 0 & 0\\
		\beta & 0 & \frac{a_2^2}{\rho_2} & 0 & 0 & 0 & u_2^1 & -\vartheta_{1}^2 & -\vartheta_{1}^3 \\
		0 & 0 & 0 & 0 & 0 & 0 & 0 & u^1 & 0  \\
		0 & 0 & 0 & 0 & 0 & 0 & 0 & 0 & u^1  \\
	\end{array}
	\right).
\end{equation*}
As previously, we have used the auxiliary variables~\eqref{eq:auxiliary_var_jacobian} to lighten the matrix. In this case, the nine eigenvalues of the matrix $\bC$ are given by
\begin{align*}
	&\lambda_1 = u^1_1-a_1,   \quad    
	 \lambda_2 = u^1_1+a_1,   \quad 
	 \lambda_{3} = u^1_2-a_2, \quad 
	 \lambda_{4} = u^1_2+a_2, \quad 
	 \lambda_{5-9} = u^1, 
\end{align*}
with $a_j$ the sound speed of phase $j$, which is defined as $a^2_j=\rho_j\left(\frac{\partial h_j}{\partial \rho_j}\right), \, j=1,2$. Since the eigenvalues are all real, to check if the system is weakly or strongly hyperbolic, it is necessary to compute the associated eigenvectors. Below, we will write the matrix that contains the right eigenvectors in columns. They are listed in the same order as the eigenvalues.
\begin{equation*}
	\bR=\left(
	\arraycolsep=2.8pt\def\arraystretch{1.5} 
	\begin{array}{ccccccccc}
		0 & 0 & 0 & 0 & \frac{\alpha_1\,\alpha_2\,\vartheta_1^3}{\eta_2} & \frac{\alpha_1\,
		\alpha_2\,\vartheta_1^2}{\eta_2} & -\frac{\alpha_1\,\alpha_2\,
		Z_{2}}{\eta_2\,\vartheta_1^1} & 0 & 0 \\[0.2cm]
		\frac{-\rho_1}{a_1} & \frac{\rho_1}{a_1} & 0 & 0 & -\frac{\zeta_1\,\rho_1\,\vartheta_1^3}{\eta_2 Z_{1}} & -\frac{\zeta_1\,\rho_1\,
		\vartheta_1^2}{\eta_2 Z_{1}} & \frac{\zeta_1\,\rho_1\,Z_{2}}{\eta_2\,\vartheta_1^1\,Z_{1}} & -\frac{\rho_1\,
		\vartheta_2^3}{Z_{1}} & -\frac{\rho_1\,\vartheta_2^2}{Z_{1}} \\[0.2cm]
		0 & 0 &	\frac{-\rho_2}{a_2} & \frac{\rho_2}{a_2} & \frac{\alpha_1\,\rho_2\,\vartheta_1^3}{\eta_2} & \frac{\alpha_1\,
		\rho_2\, \vartheta_1^2}{\eta_2} & \frac{\zeta_2\,\rho_2}{\eta_2\,\vartheta_1^1}
		& 0 & 0 \\[0.2cm]
		1 & 1 & 0 & 0 & \frac{\eta_1\,\vartheta_2^1\,\vartheta_1^3}{\eta_2\,Z_{1}} & \frac{\eta_1\,\vartheta_2^1\,\vartheta_1^2}{\eta_2\,Z_{1}} & -\frac{c_2\,\eta_1\,Z_{2}}{c_1\,\eta_2\,Z_{1}} & \frac{\vartheta_2^1 \vartheta_2^3}{Z_{1}} & \frac{\vartheta_2^1\vartheta_2^2}{Z_{1}} \\
		0 & 0 & 0 & 0 & 0 & 0 & 0 & 0 & 1 \\
		0 & 0 & 0 & 0 & 0 & 0 & 0 & 1 & 0 \\
		0 & 0 & 1 & 1 & 0 & 0 & 1 & 0 & 0 \\
		0 & 0 & 0 & 0 & 0 & 1 & 0 & 0 & 0 \\
		0 & 0 & 0 & 0 & 1 & 0 & 0 & 0 & 0 \\
	\end{array}
	\right),
\end{equation*}
where the auxiliary variables used to write the nine right eigenvectors more compactly are defined as
\begin{align*}
		&Z_{1}= a_1^2-(\vartheta_{2}^1)^2, \quad 
		Z_{2}=a_2^2-(\vartheta_{1}^1)^2,  \quad 
		\eta_1= \alpha_2(a_1^2-\alpha_1\beta), \\ 
		&\eta_2 = \alpha_1(a_2^2+\alpha_2\beta), \quad 
		\zeta_1= \alpha_2\left(\alpha_1\beta-(\vartheta_2^1)^2\right), \quad 
		\zeta_2 = \alpha_1\left(\alpha_2\beta+(\vartheta_1^1)^2\right).
\end{align*}  
Since we have obtained nine real eigenvalues with a full set of linearly independent eigenvectors, the system~\eqref{eq:GodPowell_barotropic_system} with the additional symmetrizing Godunov-Powell terms is strongly hyperbolic.  

\section{High-order ADER discontinuous Galerkin finite element scheme with \aposteriori subcell finite volume limiter}
\label{sec:dg.method}
As described in Section~\ref{sec:gov_eq} and references~\cite{Romenski2007TwoPhase,Romenski2010TwoPhase,Thein2022}, the system~\eqref{eq:cons_barotropic_system} is a hyperbolic system that can be compactly written as
\begin{equation}
\frac{\partial \bQ}{\partial t} + \nabla \cdot \bF(\bQ) + \bB{}(\bQ)\cdot \nabla \bQ = \bS(\bQ),
\label{eq:hyp_system}
\end{equation}
where $\bQ$ is the vector of conserved variables, $\bF= (\ff, \gg, \hh)$ is the  flux tensor, $\bB{}(\bQ)\cdot \nabla \bQ=\bB{1}(\bQ)\frac{\partial}{\partial x}\bQ + \bB{2}(\bQ)\frac{\partial}{\partial y}\bQ + \bB{3}(\bQ)\frac{\partial}{\partial z}\bQ$ are the non-conservative terms and $\bS$ is the vector of algebraic source terms. To solve the system~\eqref{eq:hyp_system}, we will use high-order ADER discontinuous Galerkin schemes with \aposteriori sub-cell finite volume limiter on uniform Cartesian meshes, see~\cite{Dumbser2008DG,DumbserBalsara2009PNPM,DumbserZanotti2009PNPM,Dumbser2010PNPM,Dumbser2014,Zanotti2015} for further details.

\subsection{One-Step ADER-DG Schemes}
 In the following, a description of the method is given for the two-dimensional case ($d=2$). The computational domain $\Omega = [-L_{x}/2,L_{x}/2]\times[-L_{y}/2,L_{y}/2]$ is discretized with a Cartesian grid composed of $N_{x}\times N_{y}$ cells. These cells are given by $\Omega_i = [x_i-\frac{\Delta x}{2},x_i+\frac{\Delta x}{2}]\times[y_i-\frac{\Delta x}{2},y_i+\frac{\Delta y}{2}]$, with $(x_i,y_i)$ the barycenter of the cell $\Omega_i$, and $\Delta x=\frac{L_x}{N_x}$, $\Delta y=\frac{L_y}{N_y}$ the mesh spacing in the $x-$ and $y-$directions.
Let $\bu_h(\bx,t^n)$ be the discrete solution of~\eqref{eq:hyp_system} in each spatial control volume $\Omega_i$ at time $t^n$, written in terms of tensor products of piecewise polynomials of degree $N$, and let $V_h$ be the space of tensor products of piecewise polynomials of the degree up to $N$. Then, the discrete solution $\bu_h$ can be written in terms of the basis functions, $\varphi_l(x,y), \  l\in[1,(N+1)^d]$, in every cell $\Omega_i$ as
\begin{equation}
\bu_h(\bx, t^n) = \varphi_l(\bx)\hat{\bu}_l^n,\quad \bx \in \Omega_i,
\label{eq:discrete_solution}
\end{equation}
where $\varphi_l=\varphi_l(\bx)$ are the basis functions associated with $V_h$. We take an orthogonal nodal basis $\{\varphi_l\}_{l\in\{0,\dots,(N+1)^d\}}$, generated by the tensor product $\{\varphi_{l_1} \varphi_{l_2} \varphi_{l_3} \}_{l_1,l_2,l_3\in\{0,\dots,N\}}$ where $\{\varphi_{l_l}\}_{l_l,\in\{0,\dots,N\}}$ are the Lagrange interpolation polynomials going through the $N+1$ Gauss-Legendre quadrature nodes. 
Multiplying~\eqref{eq:hyp_system} by a test function $\varphi_l\in V_h$, and integrating the equation over the space-time control volume $\Omega_i\times[t^n,t^{n+1}]$, the weak formulation can be written as
\begin{equation}
\int_{t^n}^{t^{n+1}}\!\!\!\!\int_{\Omega_i}\left(\frac{\partial \bQ}{\partial t}  + \nabla\cdot\bF(\bQ) + \bB{}(\bQ)\cdot \nabla \bQ\right)\varphi_l\,\ddx{\bx}\,\ddx{t} = \int_{t^n}^{t^{n+1}}\!\!\!\!\int_{\Omega_i} \bS(\bQ)\varphi_l\,\ddx{\bx}\,\ddx{t}.
\label{eq:weak_form}
\end{equation}
To achieve high-order in space and time, an ADER approach can be used. This methodology was put forward by Toro \textit{et al.} for the first time in \cite{Toro2001} for linear problems on Cartesian meshes, it can be implemented in both the finite volume and the discontinuous Galerkin finite element framework and is uniformly and arbitrarily high-order accurate in both space and time. In this work, an alternative version of the ADER approach is considered, avoiding the use of the Cauchy-Kovalevskaya procedure by using a local space-time discontinuous Galerkin predictor, which is based on a weak space-time formulation of the governing PDE~\eqref{eq:weak_form}.

Using~\eqref{eq:discrete_solution} and integrating the term with the time derivative by parts in time and the divergence term by parts in space, then \eqref{eq:weak_form} results in 
	\begin{align}
& \left( \int \limits_{\Omega_i} \! \varphi_k \varphi_l \, \ddx{\bx} \right) \left( \hat{\bu}_{l,i}^{n+1} - \hat{\bu}_{l,i}^{n} \, \right) +  
\int \limits_{t^n}^{t^{n+1}} \!\!\! \int \limits_{\partial\Omega_i } \! \varphi_k \left( \mathcal{G}\left(\bq_h^-, \bq_h^+ \right)  + 
\mathcal{D}\left(\bq_h^-,\bq_h^+ \right) \right) \cdot \bn \, \ddx{S}\,\ddx{t} \nonumber \\ 
&-\int \limits_{t^n}^{t^{n+1}} \!\! \int \limits_{\Omega_i } \!\! \nabla\varphi_k \cdot \bF(\bq_h) \,\ddx{\bx}\,\ddx{t}  
\int \limits_{\Omega_i} \varphi_k \bS(\bq_h) \,\ddx{\bx}\,\ddx{t},  
	\label{eq:ADER-DG_discrete}
	\end{align}
	where $\bn$ is the outward unit normal vector at the cell boundary $\partial\Omega_i$, $\bq_h$ is a local space-time predictor, which will be explained below, $\bq_h^+$ and $\bq_h^-$ are the boundary-extrapolated values of the space-time predictor from within $\Omega_i$ and its neighbor $\Omega_j$. 
Usually, $\bq_h$ presents jumps across the cell boundaries which can be resolved by the solution of a generalized Riemann problem (see~\cite{Gassner2007,Toro2009RS} for more details). In~\eqref{eq:ADER-DG_discrete}, $\mathcal{G}$ denotes the Riemann solver (numerical flux function), which depends on the left state $\bq_h^{-}$ and the right state $\bq_h^{+}$. In this case, this integral has been approximated by using the Rusanov flux, see~\cite{Rusanov1962Nonstationary}	
	\begin{equation}
	\mathcal{G}(\bq_h^-, \bq_h^+ ) \cdot \bn 
	= \frac{1}{2} \left( \bF(\bq_h^+) + \bF(\bq_h^-) \right) \cdot \bn
	- \frac{1}{2} s_{\max} \,\mathbf{I} \, \left( \bq_h^+ - \bq_h^- \right),
	\label{eq:rusanov_flux} 
	\end{equation} 
	where $s_{max} = \mathrm{max}\left(\left\vert  \lambda_k(\bq_h^+)\right\vert ,\left\vert \lambda_k(\bq_h^-) \right\vert\right)$ is the maximum wave speed at the interface.
	To deal with the jump terms in the non-conservative product, a path-conservative method is employed, following~\cite{Castro2006,Pares2006}. In this setting, a straight-line segment path is chosen 
	\begin{equation*}
		\Psi(s, \bq_h^+, \bq_h^-) =  \bq_h^-+s( \bq_h^+- \bq_h^+), \quad s\in[0,1],
	\end{equation*}
	so that the non-conservative terms reduce to the following expression
	\begin{equation*}
		\mathcal{D}(\bq_h^{-},\bq_h^{+})\cdot \bn = \frac{1}{2}\widetilde{\bB{}} \cdot (\bq_h^{+}-\bq_h^{-}), \quad \text{with} \quad \widetilde{\bB{}} = \int_0^1 \bB{}(	\Psi(s, \bq_h^+, \bq_h^-)) \cdot \bn \ \ddx{S}.
	\end{equation*}

	\subsection{Local space-time predictor}
In the following, we describe the local space-time predictor used to compute the coefficients $\hat{\bu}^n_{l,i}$ in Equation~\eqref{eq:ADER-DG_discrete}. 
We consider space-time basis functions $\theta_l$, that are obtained as the tensor product $\theta_l(\mathbf{x},t)=\varphi_{k_0}(t)\varphi_{k_1}(\mathbf{x})$, of the same previously introduced Lagrange interpolation polynomials, just that now the basis functions also depend on time. 
The predictor $\bq_h$ is written in the form
	\begin{equation}
		\bq_h(\bx,t) = \theta_l(\bx,t) \hat{\bq}_{l,i}
		\label{eq:predictor_basis},
	\end{equation}
as a weak solution to~\eqref{eq:hyp_system}. Then, using~\eqref{eq:predictor_basis} in~\eqref{eq:hyp_system}, multiplying by a space-time basis function $\theta_l$ and integrating over $\Omega_i\times[t^n,t^{n+1}]$ yields
	\begin{equation*}
	\int \limits_{t^n}^{t^{n+1}} \! \int \limits_{\Omega_i}\theta_l\left(\bU_t + \nabla \cdot \bF(\bU) + \bB{}(\bU)\cdot \nabla \bU \right) \ddx{\bx}\,\ddx{t} = \int \limits_{t^n}^{t^{n+1}} \! \int \limits_{\Omega_i}\theta_l \left(\bS(\bU)\right) \ddx{\bx}\,\ddx{t}.
	\end{equation*}
	Integrating by parts only the first term on the left-hand side and taking into account that at time $t^n$ we start from the known state $\bu_h^n$ allows us to write
	\begin{align}
	&\int \limits_{\Omega_i}\theta_l(\bx,t^{n+1}) \bq_h(\bx,t^{n+1}) \,\ddx{\bx}  
	-\int \limits_{\Omega_i}\theta_l(\bx,t^{n}) \bu_h(\bx,t^{n}) \,\ddx{\bx} 
	-\int \limits_{t^n}^{t^{n+1}} \!\! \int \limits_{\Omega_i}\frac{\partial\theta_l}{\partial t} \, \bq_h \,\ddx{\bx} \, \ddx{t} \nonumber \\
	+&\int \limits_{t^n}^{t^{n+1}} \!\! \int \limits_{\Omega_i}\theta_l\,
	\nabla \cdot \bF(\bq_h) \,\ddx{\bx} \, \ddx{t}
	+\int \limits_{t^n}^{t^{n+1}} \!\! \int \limits_{\Omega_i^{\circ} } \theta_l\bB{}(\bq_h ) \cdot \nabla \bq_h \,\ddx{\bx} \, \ddx{t} =
	\int \limits_{t^n}^{t^{n+1}} \!\! \int \limits_{\Omega_i } \theta_l \bS(\bq_h) \,\ddx{\bx} \, \ddx{t}.
	\label{eq:predictor}
	\end{align}
	Equation~\eqref{eq:predictor} is a system for the unknowns $\hat{\bq}$ of the space-time predictor $\bq_h(\bx,t)$ and can be computed in terms of the spatial degrees of freedom $\hat{\bu}_l^n$. It is solved by a fixed point iteration for which convergence was proven in~\cite{Busto2020ADERGPR}. Once the predictor is known, Equation~\eqref{eq:ADER-DG_discrete} allows to compute the polynomial coefficients $\hat{\bu}^{n+1}$ in each cell by using Gaussian quadrature for the remaining integrals.

\subsection{\textit{A posteriori} subcell finite volume limiter}
\label{subsec:TVD}

Although the numerical method presented in the previous section is a high-order method, it is linear in the sense of Godunov, which means that spurious oscillations will appear in the presence of discontinuities or shock waves. To overcome this problem, we use the \aposteriori subcell limiter for high-order fully discrete one-step ADER-DG schemes presented in~\cite{Dumbser2014,Zanotti2015}. This subcell FV limiter is based on the MOOD paradigm introduced in~\cite{Clain2011,Diot2012,Diot2013} for finite volume schemes. 

The scheme described in the previous section is run over the entire domain at each time step, and a so-called candidate solution $\bu_h^*(\bx,t^{n+1})$ is obtained. Then, one checks whether the candidate solution verifies some numerical and physical detection criteria (positivity of the densities $\rho_1$ and $\rho_2$, $\alpha_1$ with values between $0$ and $1$) and whether the discrete maximum principle (DMP), \cite{Dumbser2014}, is verified. If a cell $\Omega_i$ violates any of the above criteria, that cell is flagged as a troubled cell and for the application of the subcell finite volume limiter. The limiter is denoted as \aposteriori because it is applied after the candidate solution has been computed. 

The limiter is applied in the following way: all cells $\Omega_i$ marked as troubled are subdivided into $(2N+1)^d$ subcells, which are denoted by $\Omega_{i,j}$ where $\Omega_i=\bigcup_j \Omega_{i,j}$. The discrete solution at time $t^n$ is given by the piecewise constant cell averages, denoted by $\bar{\bu}^n_{i,j}$. They are obtained from the high-order DG polynomials $\bu_h(\bx,t^n)$ by averaging using the definition of the cell average 
	\begin{equation}
	\bar{\bu}^n_{i,j} = \frac{1}{|\Omega_{i,j}|} \int_{\Omega_{i,j}}	\bu_h(\bx,t^n) \ddx{\bx}. 
	\label{eq:subcell_averages}
	\end{equation}	
It is worth noting that subdividing a high-order DG element into $2N+1$ FV subcells per space dimension does \textit{not} reduce the time step size of the overall scheme since the CFL stability condition of explicit DG schemes scales with $1/(2N+1)$ in 1D, while the maximum Courant number of finite volume methods is unity in one space dimension.

The cell averages~\eqref{eq:subcell_averages} are evolved in time using either a second-order MUSCL-Hancock-type TVD finite-volume scheme with minmod limiter, or by making use of a third-order ADER-WENO finite-volume scheme, see~\cite{Dumbser2014}, which are also both predictor-corrector methods, and thus look almost identical to the ADER-DG scheme, except for the necessary nonlinear reconstruction step. Moreover, in this case, the test function is unity, which implies that the volume integral over the flux term disappears and the volumes computed over $\Omega_i$ are replaced by the volumes over the subcells $\Omega_{i,j}$, hence 
	\begin{align}
	\left| \Omega_{i,j} \right| \left( \bar{\bu}^{n+1}_{i,j} - \bar{\bu}^{n}_{i,j}  \right)
		+ \int_{t^n}^{t^{n+1}} \! \! \int_{\partial \Omega_{i,j}}
		\left( \mathcal{G}\left( \bq_h^-, \bq_h^+ \right)  + \mathcal{D}\left( \bq_h^-, \bq_h^+ \right) \right) \cdot \bn \, \ddx{S} \ddx{t}  \nonumber \\  
		+ \int_{t^n}^{t^{n+1}} \! \! \int_{\Omega_{i,j}^\circ}
		\left( \bB{}(\bq_h) \cdot \nabla \bq_h  \right)  \ddx{\bx} \ddx{t} 
		= \int_{t^n}^{t^{n+1}} \! \! \int_{\Omega_{i,j}}
		\bS(\bq_h, \nabla \bq_h)  \ddx{\bx} \ddx{t}\, .
	\label{eq:integral_pde}
	\end{align}
	
The limited DG polynomial $\bu'_h$ at time $t^{n+1}$ is then obtained by performing a constrained least squares reconstruction and using the averages of all the subcells of $\Omega_i$ computed using~\eqref{eq:integral_pde}. The reconstruction reads 
	\begin{equation*}
		\frac{1}{|\Omega_{i,j}|} \int_{\Omega_{i,j}} \bu'_h(\bx,t^{n+1}) \ddx{\bx}
		= \bar{\bu}^{n+1}_{i,j} \qquad \forall \Omega_{i,j} \in \Omega_i,
	\end{equation*}
with the linear constraint 
	\begin{equation}
		\int_{\Omega_{i}} \bu'_h(\bx,t^{n+1})\ddx{\bx}
		 = \displaystyle{\sum_{\Omega_{i,j} \in \Omega_i} |\Omega_{i,j}| \bar{\bu}^{n+1}_{i,j}}. 
	\label{eq:rec_constraint} 
	\end{equation}
The constraint \eqref{eq:rec_constraint} means conservation of the solution within the element $\Omega_i$. In addition to the expansion coefficients $\hat{\bu}_{i,l}^{n+1}$ of the limited DG polynomial, in all limited DG elements we also keep in memory the averages of the finite volume subcells  $\bar{\bu}^{n+1}_{i,j}$, as they serve as initial condition for the finite volume limiter of the subcell in case a cell is problematic also in the next time step, see~\cite{Dumbser2014}. More details about the \aposteriori subcell finite volume limiter can be found in~\cite{Dumbser2014,Zanotti2015,DumbserLoubere2016}.

\section{Numerical results} \label{sec:num_results}
This section is devoted to showing some test cases to illustrate the high order of accuracy of the proposed method, especially in the presence of steep gradients in the solution. First, some simulations are performed to show the experimental order of convergence (EOC) of the proposed ADER-DG method. Then, our scheme is used to solve some Riemann problems in 1D and 2D. Finally, a dambreak problem is simulated, and the results are compared with those obtained with a reduced Baer--Nunziato model. All test cases have been performed with the system~\eqref{eq:cons_barotropic_system} without using the curl-cleaning technique or the symmetrizing Godunov-Powell terms. Although the original system is only weakly hyperbolic, no stability problems have been found in the numerical simulations, contrary to what was reported in \cite{Chiocchetti2021SHTCSurfaceTension,Chiocchetti2023,Dhaouadi2023expspnsk} for other weakly hyperbolic systems with curl involutions. Moreover, in all tests, the algebraic relaxation source terms have been neglected, and the gravity $\bg$ is set to $\bm{0}$, except in the dambreak test case, where it is necessary to consider the gravity. 

\subsection{Accuracy analysis} 
This section performs a numerical convergence analysis to show the experimental order of convergence of the proposed ADER-DG method. To construct an exact solution, following~\cite{Balsara2004,Dumbser2010Force2}, an analytical, stationary, and rotationally symmetric solution of the system~\eqref{eq:cons_barotropic_system} is computed considering cylindrical coordinates ($r$, $\theta$ and $z$), with $(\velr,\velt,\velz)$ the velocity vector and $(\urelr,\urelt,\urelz)$ the relative velocity vector. The analytical solution is assumed to approach a constant state as the radial coordinate $r$ tends to infinity to be compatible with periodic boundary conditions. 

To obtain a steady analytical solution, we first write an equivalent PDE system in the radial direction. For this purpose, the pressure and the velocity relaxation are neglected, the gravity is set to zero and system~\eqref{eq:cons_barotropic_system} is rewritten in cylindrical coordinates and assuming no variations in the $z-$direction ($\partial/\partial z=0$) and considering rotational symmetry ($\partial/\partial \theta=0$). The resulting system in radial direction reads  
\begin{subequations}\label{eq:barotropic_radial}
\begin{align}
&\frac{\dalpha{1}}{\dt}+\velr \frac{\dalpha{1}}{\der r}=0,\label{eq:volume_frac_cyl}\\
&\frac{\der\alpha_1\rho_1}{\dt}+\frac{1}{r}\frac{\der (r\alpha_1\rho_1 \velr_1)}{\der r}=0,\label{eq:mass_density1_cyl}\\
&\frac{\der\alpha_2\rho_2}{\dt}+\frac{1}{r}\frac{\der (r\alpha_2\rho_2 \velr_2)}{\der r}=0,\label{eq:mass_density2_cyl}\\
&\frac{\der\rho \velr}{\dt}+\frac{\der(\alpha_1\rho_1 (\velr_1)^2 + \alpha_2\rho_2 (\velr_2)^2 + p )}{\der r}+\alpha_1\rho_1\frac{ (\velr_1)^2-(\velt_1)^2}{ r}\nonumber\\
&\hspace{5.9cm} +\alpha_2\rho_2\frac{  (\velr_2)^2-(\velt_2)^2 }{ r}=0,\label{eq:momentum1_cyl}\\
&\frac{\der\rho \velt}{\dt}+\frac{\der(\alpha_1\rho_1 \velr_1\velt_1 + \alpha_2\rho_2 \velr_2\velt_2 )}{\der r}+2\frac{\alpha_1\rho_1 \velr_1 \velt_1 + \alpha_2\rho_2 \velr_2\velt_2 }{ r}=0,\label{eq:momentum2_cyl}\\
&\frac{\der \urelr}{\dt}+\frac{\der}{\der r} \left( \frac{1}{2}(u_1^l)^2-\frac{1}{2}(u_2^l)^2 + h_1-h_2\right) =0,\label{eq:rel_velocity1_cyl}\\
&\frac{\der \urelt}{\dt}=0,\label{eq:rel_velocity2_cyl}
\end{align}
\end{subequations}
where the constraint $\nabla \times \burel=0$ has been used in the last two equations. Since we are looking for a vortex-type solution, the radial velocities vanish, that is, we set $\velr=\velr_1=\velr_2=\urelr=0$. Also, we are interested in a stationary solution, hence $\partial_t=0$. With these assumptions, the system~\eqref{eq:barotropic_radial} reduces to
\begin{subequations}\label{eq:radial_stationary}
	\begin{align}
		\alpha_1\rho_1\frac{(\velt_1)^2}{ r}+\alpha_2\rho_2\frac{  (\velt_2)^2 }{ r}-\frac{\der p}{\der r}&=0,\label{eq:momentum1_cyl_stationary}\\
		\frac{\der}{\der r} \left( \frac{1}{2}(\velt)^2-\frac{1}{2}(\velt)^2 + h_1-h_2\right) &=0.\label{eq:rel_velocity1_cyl_stationary}
	\end{align}
\end{subequations}
With some simple algebra in~\eqref{eq:radial_stationary} we get
\begin{subequations}\label{eq:velocity_vortex}
	\begin{align}
		&(\velt_1)^2=\frac{r}{\rho}\frac{\der p}{\der r}+2\frac{\alpha_2\rho_2}{\rho}(k-h_1+h_2),\\
		&(\velt_2)^2=\frac{r}{\rho}\frac{\der p}{\der r}-2\frac{\alpha_1\rho_1}{\rho}(k-h_1+h_2), 
	\end{align}
\end{subequations}
where $k$ is a constant. Now, we \textit{prescribe} radial profiles for $\alpha_1$ and $p_i$ as 
\begin{align*}
	\alpha_1 = \frac{1}{3}+\frac{e^{-\frac{ r^2}{2}}}{2\sqrt{2\pi}}, \qquad 
	p_i = 1-\frac{e^{1-r^2}}{4},
\end{align*}
and then the densities $\rho_i$ and the velocities $u_{i\theta}$ result as 
\begin{align*}
	\rho_i = \left(1 - \frac{e^{1 - r^2}}{4} \right)^{5/7}, \qquad 
	\bu_{i\theta} = 2^{3/14}\sqrt{ \frac{e^{1 - r^2}r^2}{\left(4 - e^{1 - r^2}\right)^{5/7}}}.
\end{align*}
With this, we have computed an exact stationary and rotationally symmetric solution of the PDE system~\eqref{eq:barotropic_radial} that approaches a constant state as $r\rightarrow \infty$ to be compatible with periodic boundary conditions. Then, using the principle of Galilean invariance, we can make the test unsteady if we add a uniform velocity field to this solution. After one advection period through a periodic computational domain, the exact solution will be given by the initial condition and we can perform the convergence test.
To analyze the convergence order, we compute the solution with the proposed method, using different orders for the DG scheme, and compare it with the exact solution derived above. The computational domain is $\Omega=[-10,10]^2$ with final simulation time $t=1$ and periodic boundary conditions everywhere. Different polynomial approximation degrees have been considered for the DG scheme. The $L^{2}$ errors and the corresponding numerical convergence rates for $N=2,3,4,5$, are given in Table~\ref{tab:CT_v0}, showing the expected order of convergence. 

For analyzing the convergence order using the unsteady solution, we make use of the principle of Galilean invariance of Newtonian mechanics. We add a constant uniform velocity field $\bar{\vel}_1=\bar{\vel}_2=4$ to both phases. The computational domain is the same as before $\Omega=[-10,10]\times[-10,10]$ with final simulation time $t=5$ and periodic boundary conditions everywhere. The $L^{2}$ errors and the corresponding convergence rates for the different degrees $N=2,3,4,5$,  are given in Table~\ref{tab:CT_v4}, finding the expected convergence order $N+1$ of our high-order ADER-DG schemes.
\begin{table}[!h]
    \centering
    \scriptsize
    \caption{Numerical convergence results for high-order DG schemes of polynomial approximation $N=2,3,4,5$ with a uniform Cartesian mesh of $N_x\times N_y$ elements. The $L^{2}$ error norms and the corresponding orders of convergence of the variables $\alpha_1$, $\rho_1$, $\rho_2$, $\vel_1$, and $\vel_2$ are computed at time $t=1$.}
    \begin{tabular}{cccccc}
    \toprule
    \multicolumn{6}{c}{$N=2$}\\[-0.05cm] \midrule
	$N_x=N_y$ & $L^{2}_{\Omega}\left(\alpha_1\right)$ & $L^{2}_{\Omega}\left(\rho_1\right)$ & $L^{2}_{\Omega}\left(\rho_2\right)$ & $L^{2}_{\Omega}\left(u_1\right)$ & $L^{2}_{\Omega}\left(u_2\right)$   \\[-0.05cm] \midrule
        $16$ & $4.2737\cdot 10^{-3}$ & $2.5907\cdot 10^{-2}$ & $2.6117\cdot 10^{-2}$ & $6.7706\cdot 10^{-2}$  & $6.4058\cdot 10^{-2}$   \\ 
        $32$ & $1.0776\cdot 10^{-3}$ & $4.7315\cdot 10^{-3}$ & $4.4746\cdot 10^{-3}$ & $1.0146\cdot 10^{-2}$ & $9.2257\cdot 10^{-3}$  \\ 
        $64$ & $2.1380\cdot 10^{-4}$ & $6.7212\cdot 10^{-4}$ & $6.2405\cdot 10^{-4}$ & $1.6479\cdot 10^{-3}$ & $1.5250\cdot 10^{-3}$   \\ 
        $128$ & $3.3522\cdot 10^{-5}$ & $9.5237\cdot 10^{-5}$ & $8.7193\cdot 10^{-5}$ & $2.7709\cdot 10^{-4}$  & $2.6116\cdot 10^{-4}$  \\ 
        $256$ & $4.7854\cdot 10^{-6}$ & $1.2987\cdot 10^{-5}$ & $1.1645\cdot 10^{-5}$ & $4.3836\cdot 10^{-5}$ & $4.1955\cdot 10^{-5}$  \\[-0.05cm] \midrule
   & $\mathcal{O}(\alpha_1)$ & $\mathcal{O}(\rho_1)$ & $\mathcal{O}(\rho_2)$ & $\mathcal{O}(u_1)$ & $\mathcal{O}(u_2)$ \\[-0.05cm] \midrule
         & $1.99$ & $2.45$ & $2.55$ & $2.74$ & $2.80$   \\ 
         & $2.33$ & $2.82$ & $2.84$ & $2.62$ & $2.60$   \\ 
         & $2.67$ & $2.82$ & $2.84$ & $2.57$ & $2.55$   \\ 
         & $2.81$ & $2.87$ & $2.90$ & $2.66$ & $2.64$   \\[-0.05cm]
		\toprule
		\multicolumn{6}{c}{$N=3$}\\[-0.05cm] \midrule
		$N_x=N_y$ & $L^{2}_{\Omega}\left(\alpha_1\right)$ & $L^{2}_{\Omega}\left(\rho_1\right)$ & $L^{2}_{\Omega}\left(\rho_2\right)$ & $L^{2}_{\Omega}\left(u_1\right)$ &  $L^{2}_{\Omega}\left(u_2\right)$ \\[-0.05cm] \midrule
		$16$ & $8.2842\cdot 10^{-4}$ & $6.2523\cdot 10^{-3}$ & $4.2100\cdot 10^{-3}$ & $1.3507\cdot 10^{-2}$ & $1.1535\cdot 10^{-2}$ \\ 
		$32$ & $3.4241\cdot 10^{-5}$ & $2.5841\cdot 10^{-4}$ & $2.6524\cdot 10^{-4}$ & $1.4769\cdot 10^{-3}$ & $1.4008\cdot 10^{-3}$ \\ 
		$64$ & $1.4470\cdot 10^{-6}$ & $1.3027\cdot 10^{-5}$ & $1.0082\cdot 10^{-5}$ & $6.9613\cdot 10^{-5}$ & $5.5416\cdot 10^{-5}$ \\ 
		$96$ & $2.2214\cdot 10^{-7}$ & $2.3498\cdot 10^{-6}$ & $1.6163\cdot 10^{-6}$ & $1.1279\cdot 10^{-5}$ & $8.2417\cdot 10^{-6}$ \\ 
		$128$ & $5.9489\cdot 10^{-8}$ & $7.6189\cdot 10^{-7}$ & $5.4790\cdot 10^{-7}$ & $3.2599\cdot 10^{-6}$ & $2.2986\cdot 10^{-6}$\\[-0.05cm] \midrule
		& $\mathcal{O}(\alpha_1)$ & $\mathcal{O}(\rho_1)$ & $\mathcal{O}(\rho_2)$ & $\mathcal{O}(u_1)$ & $\mathcal{O}(u_2)$ \\[-0.05cm] \midrule
		~ & $4.60$ & $4.60$ & $3.99$ & $3.19$ & $3.04$ \\ 
		~ & $4.56$ & $4.31$ & $4.72$ & $4.41$ & $4.66$ \\ 
		~ & $4.62$ & $4.22$ & $4.51$ & $4.49$ & $4.70$ \\ 
		~ & $4.58$ & $3.92$ & $3.76$ & $4.31$ & $4.44$ \\[-0.05cm] 
		\toprule
		\multicolumn{6}{c}{$N=4$}\\[-0.05cm] \midrule
		$N_x=N_y$ & $L^{2}_{\Omega}\left(\alpha_1\right)$ & $L^{2}_{\Omega}\left(\rho_1\right)$ & $L^{2}_{\Omega}\left(\rho_2\right)$ & $L^{2}_{\Omega}\left(u_1\right)$ &  $L^{2}_{\Omega}\left(u_2\right)$ \\[-0.05cm] \midrule
		$16$ & $1.0984\cdot 10^{-4}$ & $7.6397\cdot 10^{-4}$ & $8.7616\cdot 10^{-4}$ & $4.6037\cdot 10^{-3}$ & $4.3124\cdot 10^{-3}$ \\ 
		$24$ & $1.7440\cdot 10^{-5}$ & $1.4231\cdot 10^{-4}$ & $1.1489\cdot 10^{-4}$ & $7.2792\cdot 10^{-4}$ & $5.8106\cdot 10^{-4}$ \\ 
		$32$ & $5.1768\cdot 10^{-6}$ & $3.9079\cdot 10^{-5}$ & $2.8239\cdot 10^{-5}$ & $1.8850\cdot 10^{-4}$ & $1.3895\cdot 10^{-4}$ \\ 
		$48$ & $8.9754\cdot 10^{-7}$ & $5.8266\cdot 10^{-6}$ & $4.2059\cdot 10^{-6}$ & $2.4454\cdot 10^{-5}$ & $1.8810\cdot 10^{-5}$ \\ 
		$64$ & $2.4317\cdot 10^{-7}$ & $1.4682\cdot 10^{-6}$ & $1.0720\cdot 10^{-6}$ & $5.7687\cdot 10^{-6}$ & $4.7316\cdot 10^{-6}$\\[-0.05cm] \midrule
		& $\mathcal{O}(\alpha_1)$ & $\mathcal{O}(\rho_1)$ & $\mathcal{O}(\rho_2)$ & $\mathcal{O}(u_1)$ & $\mathcal{O}(u_2)$ \\[-0.05cm] \midrule
		& $4.54$ & $4.14$ & $5.01$ & $4.55$ & $4.94$ \\ 
		& $4.22$ & $4.49$ & $4.88$ & $4.70$ & $4.97$ \\ 
		& $4.32$ & $4.69$ & $4.70$ & $5.04$ & $4.93$ \\ 
		& $4.54$ & $4.79$ & $4.75$ & $5.02$ & $4.80$ \\[-0.05cm] 
		\toprule
		\multicolumn{6}{c}{$N=5$}\\[-0.05cm] \midrule
		$N_x=N_y$ & $L^{2}_{\Omega}\left(\alpha_1\right)$ & $L^{2}_{\Omega}\left(\rho_1\right)$ & $L^{2}_{\Omega}\left(\rho_2\right)$ & $L^{2}_{\Omega}\left(u_1\right)$ &  $L^{2}_{\Omega}\left(u_2\right)$ \\[-0.05cm] \midrule
		$10$ & $3.1044\cdot 10^{-4}$ & $1.8184\cdot 10^{-3}$ & $2.2534\cdot 10^{-3}$ & $1.0731\cdot 10^{-2}$ & $9.9726\cdot 10^{-3}$ \\ 
		$20$ & $9.4868\cdot 10^{-6}$ & $6.9885\cdot 10^{-5}$ & $3.4620\cdot 10^{-5}$ & $2.6270\cdot 10^{-4}$ & $1.8766\cdot 10^{-4}$ \\ 
		$30$ & $1.3067\cdot 10^{-6}$ & $7.6256\cdot 10^{-6}$ & $6.9939\cdot 10^{-6}$ & $5.2178\cdot 10^{-5}$ & $5.0666\cdot 10^{-5}$ \\ 
		$40$ & $2.6362\cdot 10^{-7}$ & $1.4736\cdot 10^{-6}$ & $1.5563\cdot 10^{-6}$ & $1.1932\cdot 10^{-5}$ & $1.1776\cdot 10^{-5}$ \\
		$50$ & $5.6596\cdot 10^{-8}$ & $3.6572\cdot 10^{-7}$ & $3.8902\cdot 10^{-7}$ & $3.0786\cdot 10^{-6}$ & $3.0067\cdot 10^{-6}$ \\[-0.05cm] \midrule
		& $\mathcal{O}(\alpha_1)$ & $\mathcal{O}(\rho_1)$ & $\mathcal{O}(\rho_2)$ & $\mathcal{O}(u_1)$ & $\mathcal{O}(u_2)$ \\[-0.05cm] \midrule
		& $5.03$ & $4.70$ & $6.02$ & $5.35$ & $5.73$ \\ 
		& $4.89$ & $5.46$ & $3.94$ & $3.99$ & $3.23$ \\ 
		& $5.56$ & $5.71$ & $5.22$ & $5.13$ & $5.07$ \\ 
		& $6.89$ & $6.25$ & $6.21$ & $6.07$ & $6.12$ \\[-0.05cm] \bottomrule
	\end{tabular}
	\label{tab:CT_v0}
\end{table}

\begin{table}[!ht]
    \centering
    \scriptsize
    \caption{Numerical convergence rates for DG schemes of polynomial approximation $N=2,3,4,5$ with a uniform Cartesian mesh of $N_x\times N_y$ elements in the unsteady case. The $L^{2}$ error norms and the convergence orders of the variables $\alpha_1$, $\rho_1$, $\rho_2$, $\vel_1$ and $\vel_2$, are computed at time $t=5$, with $\bar{\vel}_1=\bar{\vel}_2=4.$}
    \begin{tabular}{cccccc}
    \toprule
    \multicolumn{6}{c}{$N=2$}\\[-0.05cm] \midrule
	$N_x=N_y$ & $L^{2}_{\Omega}\left(\alpha_1\right)$ & $L^{2}_{\Omega}\left(\rho_1\right)$ & $L^{2}_{\Omega}\left(\rho_2\right)$ & $L^{2}_{\Omega}\left(u_1\right)$ &  $L^{2}_{\Omega}\left(u_2\right)$\\[-0.05cm] \midrule
        $16$ & $2.5706\cdot 10^{-2}$ & $1.2916\cdot 10^{-1}$ & $1.2627\cdot 10^{-1}$ & $2.6127\cdot 10^{-1}$ & $2.5967\cdot 10^{-1}$ \\ 
        $32$ & $2.8324\cdot 10^{-3}$ & $1.4294\cdot 10^{-2}$ & $1.4025\cdot 10^{-2}$ & $3.3462\cdot 10^{-2}$ & $3.2776\cdot 10^{-2}$ \\ 
        $64$ & $3.1637\cdot 10^{-4}$ & $1.6697\cdot 10^{-3}$ & $1.6111\cdot 10^{-3}$ & $3.7392\cdot 10^{-3}$ & $3.6964\cdot 10^{-3}$ \\
        $128$ & $3.5640\cdot 10^{-5}$ & $1.9792\cdot 10^{-4}$ & $1.8772\cdot 10^{-4}$ & $4.1824\cdot 10^{-4}$ & $4.1459\cdot 10^{-4}$ \\ 
        $256$ & $4.3155\cdot 10^{-6}$ & $2.4473\cdot 10^{-5}$ & $2.3029\cdot 10^{-5}$ & $5.0492\cdot 10^{-5}$ & $5.0066\cdot 10^{-5}$ \\[-0.05cm] \midrule 
   & $\mathcal{O}(\alpha_1)$ & $\mathcal{O}(\rho_1)$ & $\mathcal{O}(\rho_2)$ & $\mathcal{O}(u_1)$ & $\mathcal{O}(u_2)$ \\[-0.05cm] \midrule
         & $3.18$ & $3.18$ & $3.17$ & $2.96$ & $2.99$ \\ 
         & $3.16$ & $3.10$ & $3.12$ & $3.16$ & $3.15$ \\ 
         & $3.15$ & $3.08$ & $3.10$ & $3.16$ & $3.16$ \\ 
         & $3.05$ & $3.02$ & $3.03$ & $3.05$ & $3.05$ \\[-0.05cm]
         \toprule
         \multicolumn{6}{c}{$N=3$}\\[-0.05cm] \midrule
	$N_x=N_y$ & $L^{2}_{\Omega}\left(\alpha_1\right)$ & $L^{2}_{\Omega}\left(\rho_1\right)$ & $L^{2}_{\Omega}\left(\rho_2\right)$ & $L^{2}_{\Omega}\left(u_1\right)$ & $L^{2}_{\Omega}\left(u_2\right)$ \\[-0.05cm] \midrule
        $16$ & $2.5156\cdot 10^{-3}$ & $1.5368\cdot 10^{-2}$ & $1.4190\cdot 10^{-2}$ & $4.0072\cdot 10^{-2}$ & $3.8626\cdot 10^{-2}$ \\
        $32$ & $1.3996\cdot 10^{-4}$ & $1.6489\cdot 10^{-3}$ & $1.5580\cdot 10^{-3}$ & $4.4715\cdot 10^{-3}$ & $4.4706\cdot 10^{-3}$ \\
        $64$ & $9.2451\cdot 10^{-6}$ & $7.2532\cdot 10^{-5}$ & $6.7138\cdot 10^{-5}$ & $2.0963\cdot 10^{-4}$ & $2.0941\cdot 10^{-4}$ \\
        $96$ & $1.8419\cdot 10^{-6}$ & $1.0446\cdot 10^{-5}$ & $9.2220\cdot 10^{-6}$ & $2.8028\cdot 10^{-5}$ & $2.7860\cdot 10^{-5}$ \\ 
        $128$ & $5.8446\cdot 10^{-7}$ & $2.9024\cdot 10^{-6}$ & $2.4794\cdot 10^{-6}$ & $7.2396\cdot 10^{-6}$ & $7.1617\cdot 10^{-6}$ \\[-0.05cm] \midrule  
   & $\mathcal{O}(\alpha_1)$ & $\mathcal{O}(\rho_1)$ & $\mathcal{O}(\rho_2)$ & $\mathcal{O}(u_1)$ &  $\mathcal{O}(u_2)$ \\[-0.05cm] \midrule
         & $4.17$ & $3.22$ & $3.19$ & $3.16$ & $3.11$ \\ 
         & $3.92$ & $4.51$ & $4.54$ & $4.41$ & $4.42$ \\ 
         & $3.98$ & $4.78$ & $4.90$ & $4.96$ & $4.97$ \\ 
         & $3.99$ & $4.45$ & $4.57$ & $4.71$ & $4.72$ \\[-0.05cm]
         \toprule
         \multicolumn{6}{c}{$N=4$} \\[-0.05cm] \midrule
	$N_x=N_y$ & $L^{2}_{\Omega}\left(\alpha_1\right)$ & $L^{2}_{\Omega}\left(\rho_1\right)$ & $L^{2}_{\Omega}\left(\rho_2\right)$ & $L^{2}_{\Omega}\left(u_1\right)$ &  $L^{2}_{\Omega}\left(u_2\right)$ \\[-0.05cm] \midrule
        $16$ & $7.2439\cdot 10^{-4}$ & $6.4964\cdot 10^{-3}$ & $6.2200\cdot 10^{-3}$ & $1.7050\cdot 10^{-2}$ & $1.6952\cdot 10^{-2}$ \\ 
        $24$ & $7.8048\cdot 10^{-5}$ & $9.1442\cdot 10^{-4}$ & $8.7219\cdot 10^{-4}$ & $2.6937\cdot 10^{-3}$ & $2.6919\cdot 10^{-3}$ \\ 
        $32$ & $1.4753\cdot 10^{-5}$ & $1.6634\cdot 10^{-4}$ & $1.5429\cdot 10^{-4}$ & $5.4339\cdot 10^{-4}$ & $5.4258\cdot 10^{-4}$ \\ 
        $48$ & $1.6068\cdot 10^{-6}$ & $1.2356\cdot 10^{-5}$ & $9.8176\cdot 10^{-6}$ & $4.0228\cdot 10^{-5}$ & $4.0030\cdot 10^{-5}$ \\ 
        $64$ & $3.8068\cdot 10^{-7}$ & $2.4863\cdot 10^{-6}$ & $1.7373\cdot 10^{-6}$ & $6.2721\cdot 10^{-6}$ & $6.2039\cdot 10^{-6}$ \\[-0.05cm] \midrule 
   & $\mathcal{O}(\alpha_1)$ & $\mathcal{O}(\rho_1)$ & $\mathcal{O}(\rho_2)$ & $\mathcal{O}(u_1)$ &  $\mathcal{O}(u_2)$\\[-0.05cm] \midrule
         & $5.49$ & $4.84$ & $4.85$ & $4.55$ & $4.54$ \\ 
         & $5.79$ & $5.92$ & $6.02$ & $5.56$ & $5.57$ \\ 
         & $5.47$ & $6.41$ & $6.79$ & $6.42$ & $6.43$ \\ 
         & $5.01$ & $5.57$ & $6.02$ & $6.46$ & $6.48$ \\[-0.05cm]
         \toprule
         \multicolumn{6}{c}{$N=5$} \\[-0.05cm] \midrule
	$N_x=N_y$ & $L^{2}_{\Omega}\left(\alpha_1\right)$ & $L^{2}_{\Omega}\left(\rho_1\right)$ & $L^{2}_{\Omega}\left(\rho_2\right)$ & $L^{2}_{\Omega}\left(u_1\right)$ &  $L^{2}_{\Omega}\left(u_2\right)$  \\[-0.05cm] \midrule
		$10$ & $2.0532\cdot 10^{-3}$ & $1.4798\cdot 10^{-2}$ & $1.3807\cdot 10^{-2}$ & $3.5831\cdot 10^{-2}$ & $3.5090\cdot 10^{-2}$ \\ 
        $20$ & $2.9701\cdot 10^{-5}$ & $3.4169\cdot 10^{-4}$ & $3.1938\cdot 10^{-4}$ & $1.1606\cdot 10^{-3}$ & $1.1572\cdot 10^{-3}$ \\ 
        $30$ & $2.2845\cdot 10^{-6}$ & $1.8935\cdot 10^{-5}$ & $1.6743\cdot 10^{-5}$ & $8.3094\cdot 10^{-5}$ & $8.2915\cdot 10^{-5}$ \\ 
        $40$ & $3.8493\cdot 10^{-7}$ & $3.1932\cdot 10^{-6}$ & $2.7274\cdot 10^{-6}$ & $1.1770\cdot 10^{-5}$ & $1.1754\cdot 10^{-5}$ \\ 
        $50$ & $9.4828\cdot 10^{-8}$ & $9.6928\cdot 10^{-7}$ & $8.3496\cdot 10^{-7}$ & $3.1726\cdot 10^{-6}$ & $3.1671\cdot 10^{-6}$ \\[-0.05cm] \midrule  
   & $\mathcal{O}(\alpha_1)$ & $\mathcal{O}(\rho_1)$ & $\mathcal{O}(\rho_2)$ & $\mathcal{O}(u_1)$ & $\mathcal{O}(u_2)$ \\[-0.05cm] \midrule
         & $6.11$ & $5.44$ & $5.43$ & $4.95$ & $4.92$ \\ 
         & $6.33$ & $7.13$ & $7.27$ & $6.50$ & $6.50$ \\ 
         & $6.19$ & $6.19$ & $6.31$ & $6.79$ & $6.79$ \\ 
         & $6.28$ & $5.34$ & $5.30$ & $5.88$ & $5.88$ \\[-0.05cm] \bottomrule 
    \end{tabular}
    \label{tab:CT_v4}
\end{table} 

\subsection{1D Riemann problems} \label{sec:RP1D} 

\begin{figure}
	\begin{center}
		\begin{tabular}{cc} 
			\includegraphics[trim= 5 5 5 5,clip,width=0.35\textwidth]{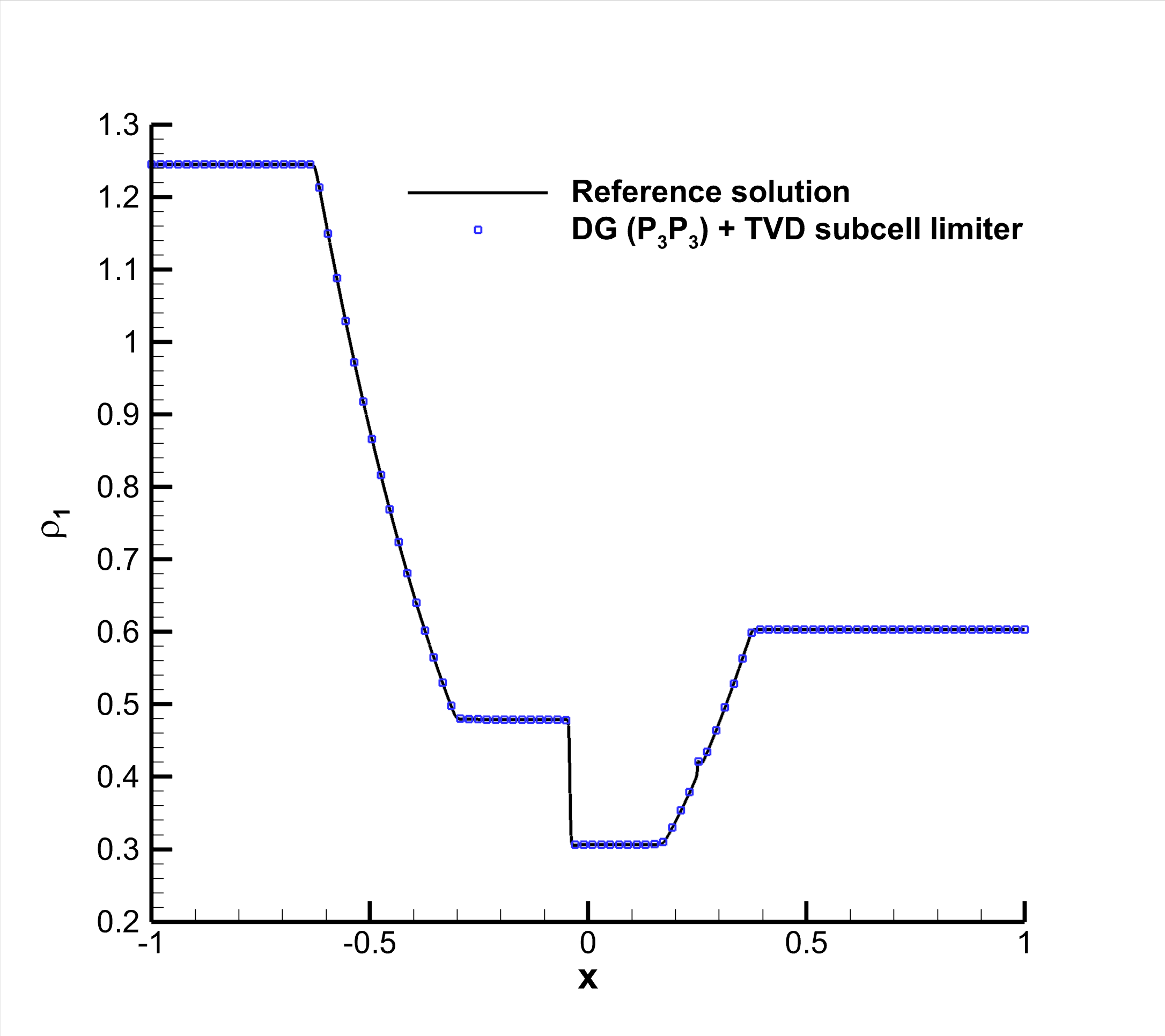}   & 
			\includegraphics[trim= 5 5 5 5,clip,width=0.35\textwidth]{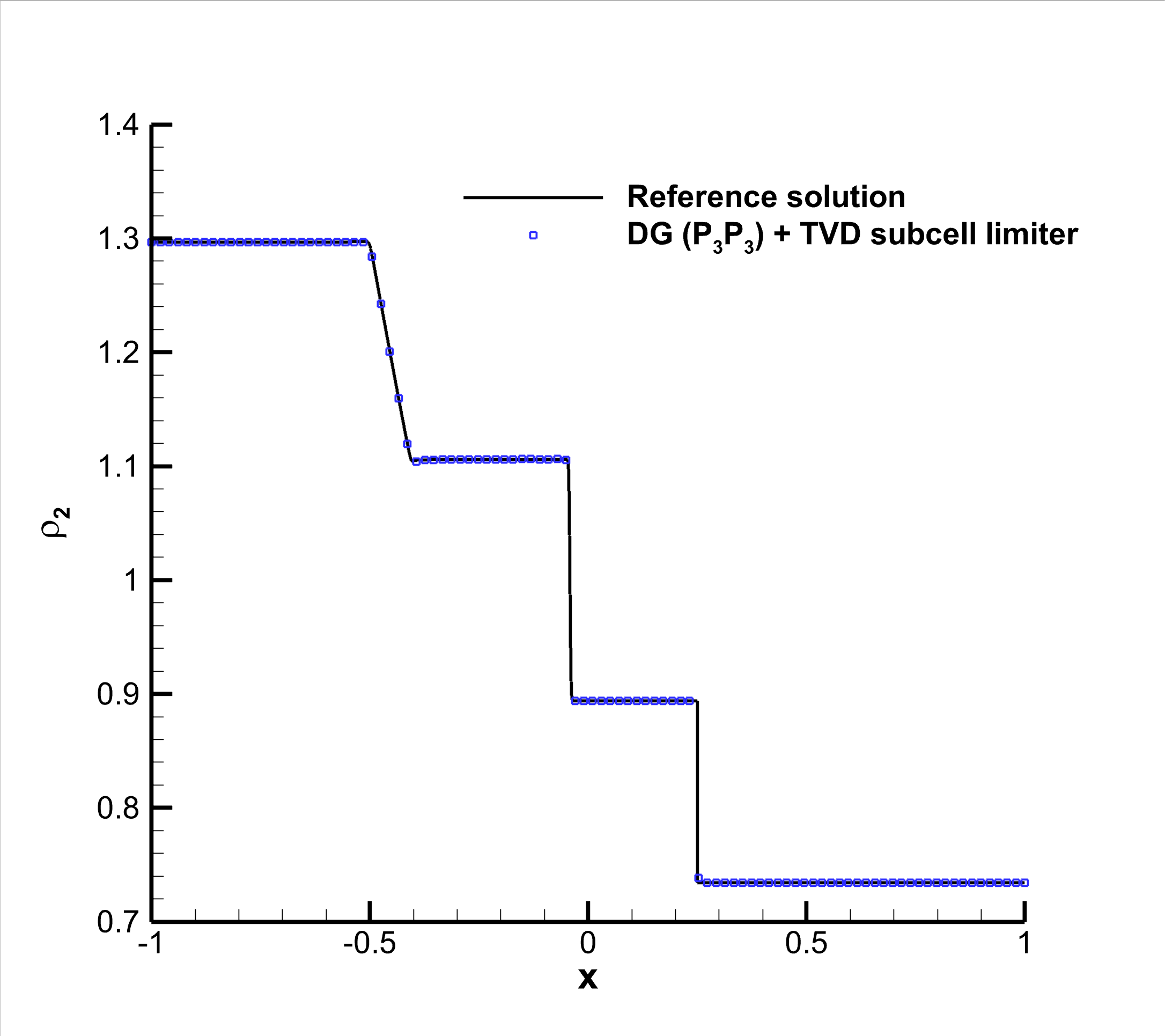}   \\[-0.2cm] 
			\includegraphics[trim= 5 5 5 5,clip,width=0.35\textwidth]{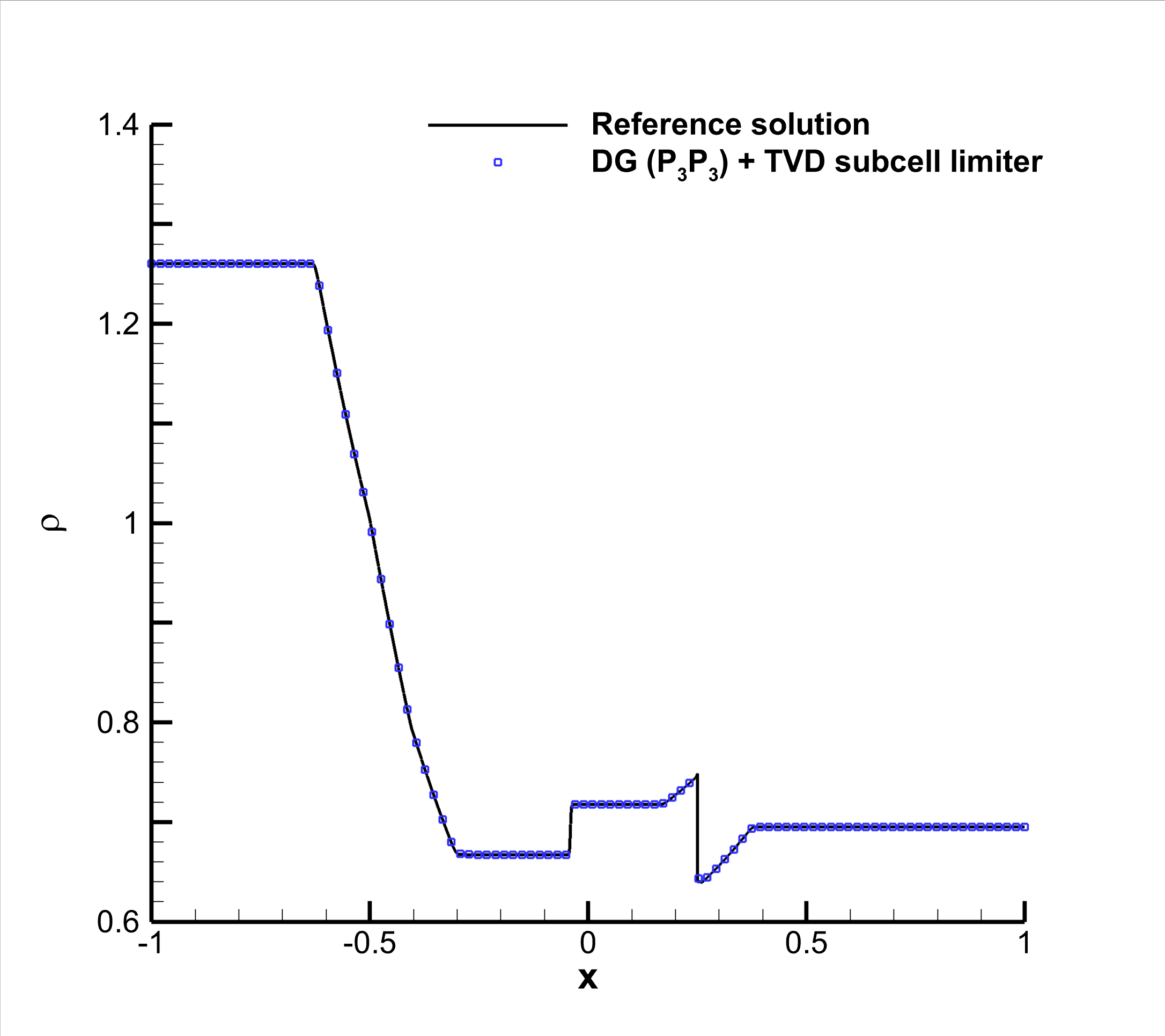}    &   
			\includegraphics[trim= 5 5 5 5,clip,width=0.35\textwidth]{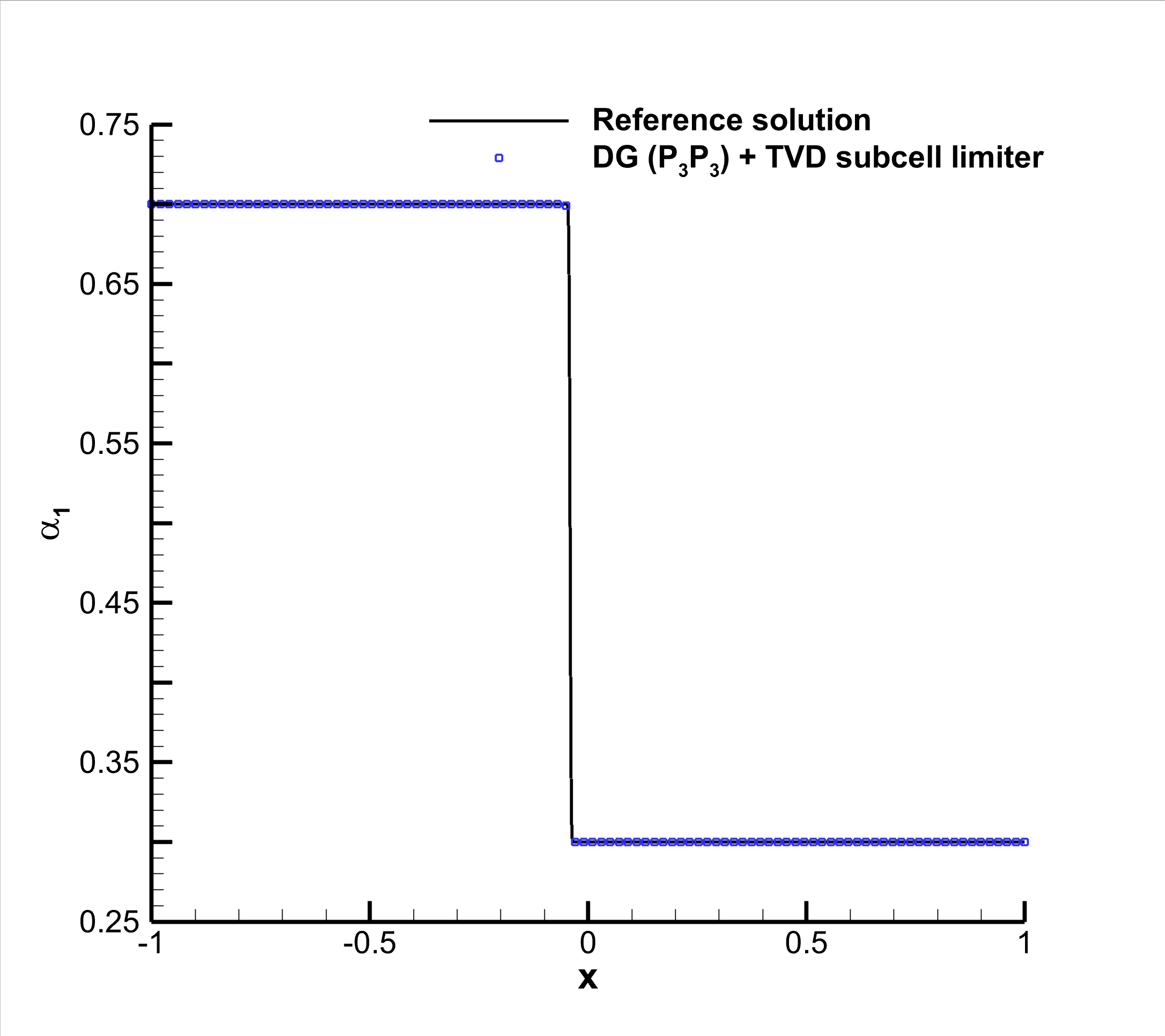}  \\[-0.2cm] 
			\includegraphics[trim= 5 5 5 5,clip,width=0.35\textwidth]{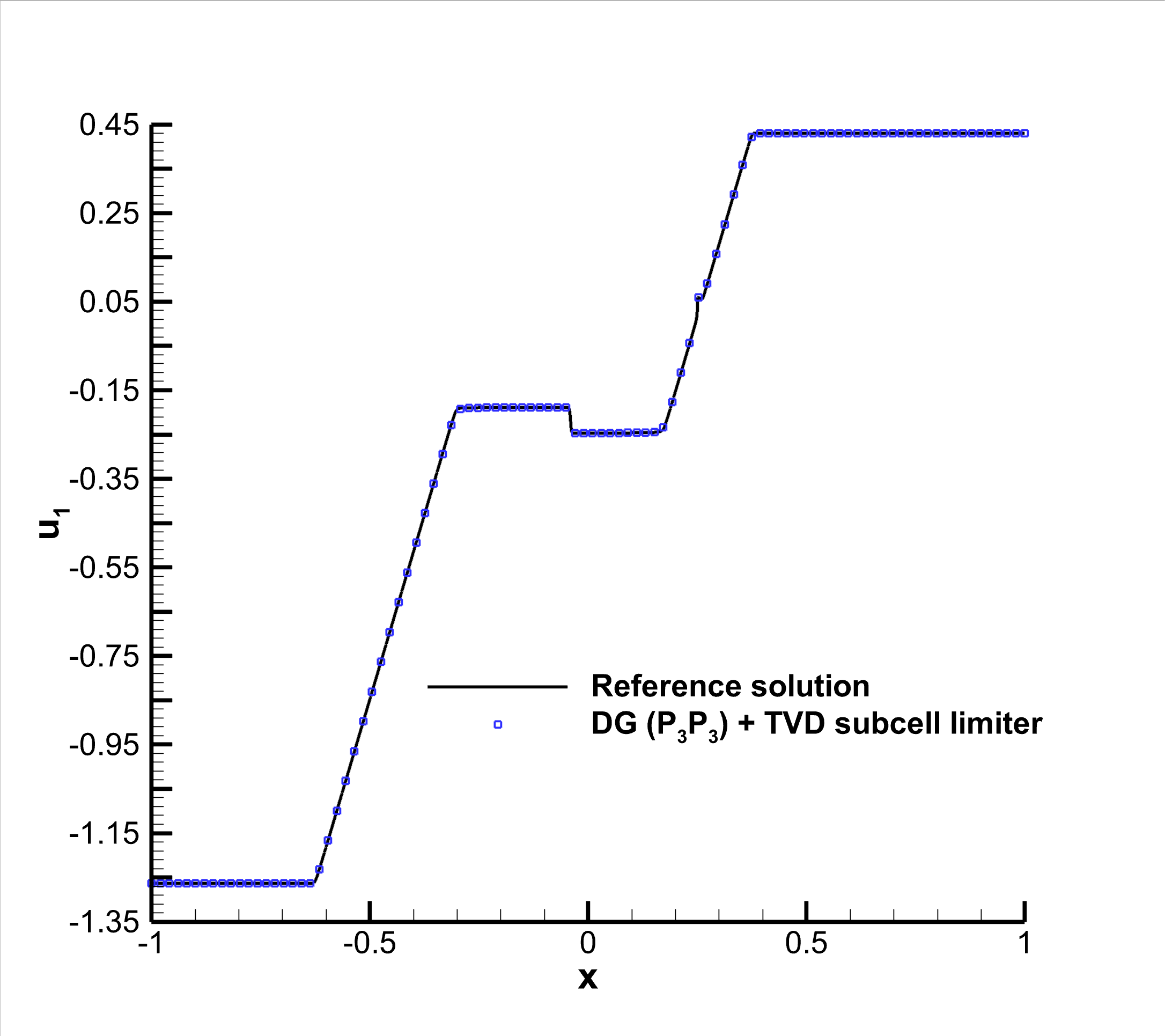}     & 
			\includegraphics[trim= 5 5 5 5,clip,width=0.35\textwidth]{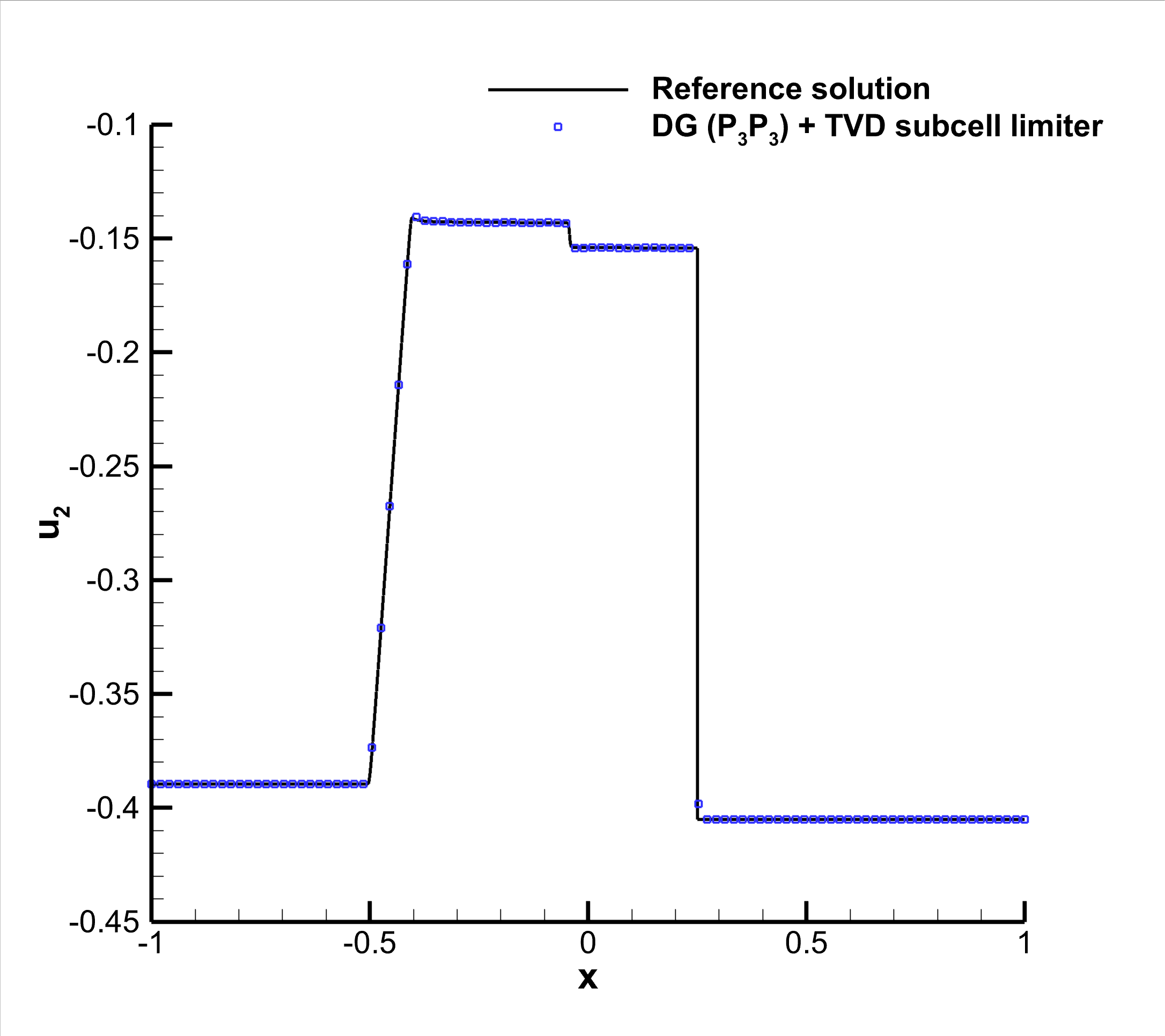}     \\[-0.2cm] 
			\includegraphics[trim= 5 5 5 5,clip,width=0.35\textwidth]{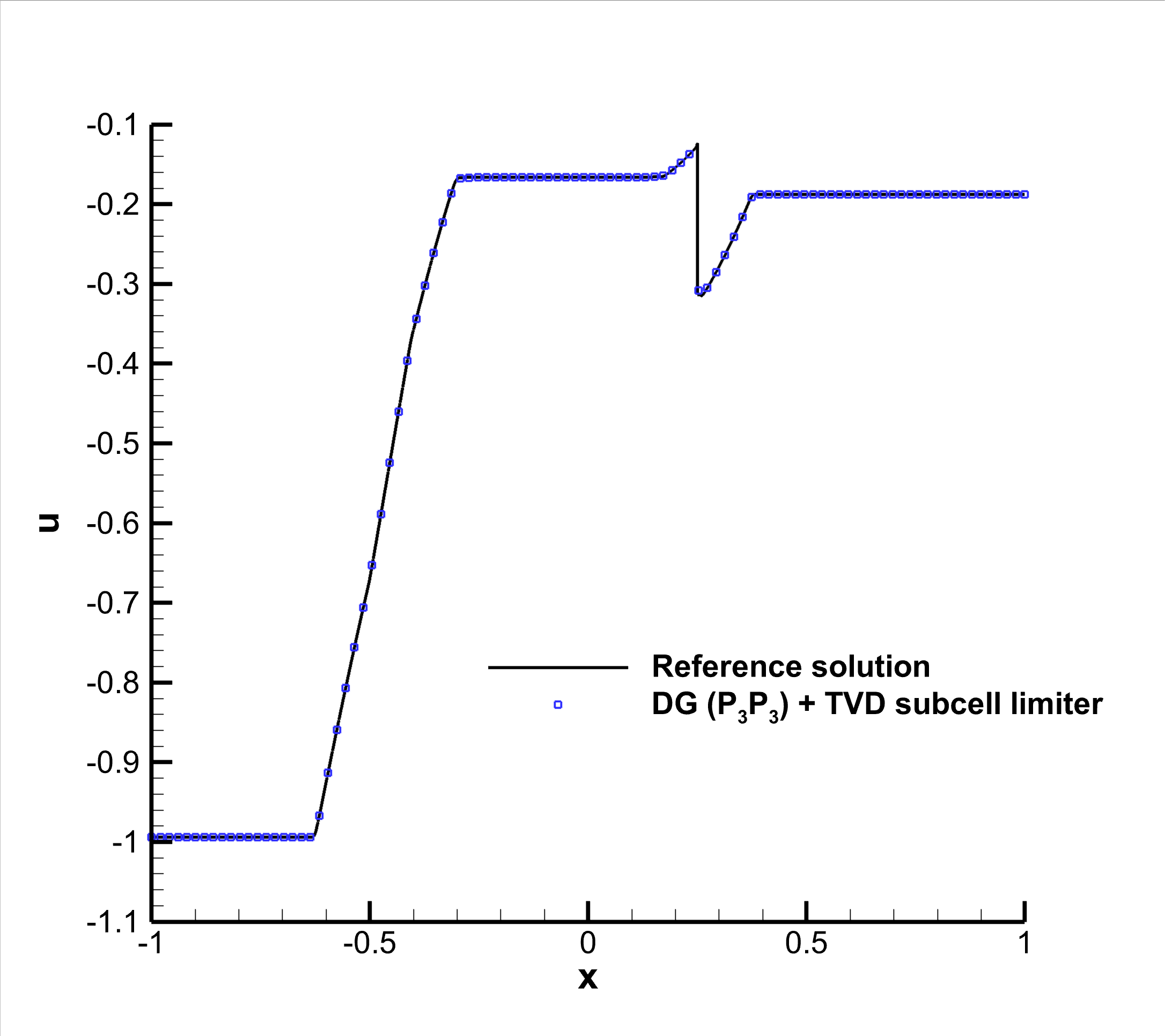}      &
			\includegraphics[trim= 5 5 5 5,clip,width=0.35\textwidth]{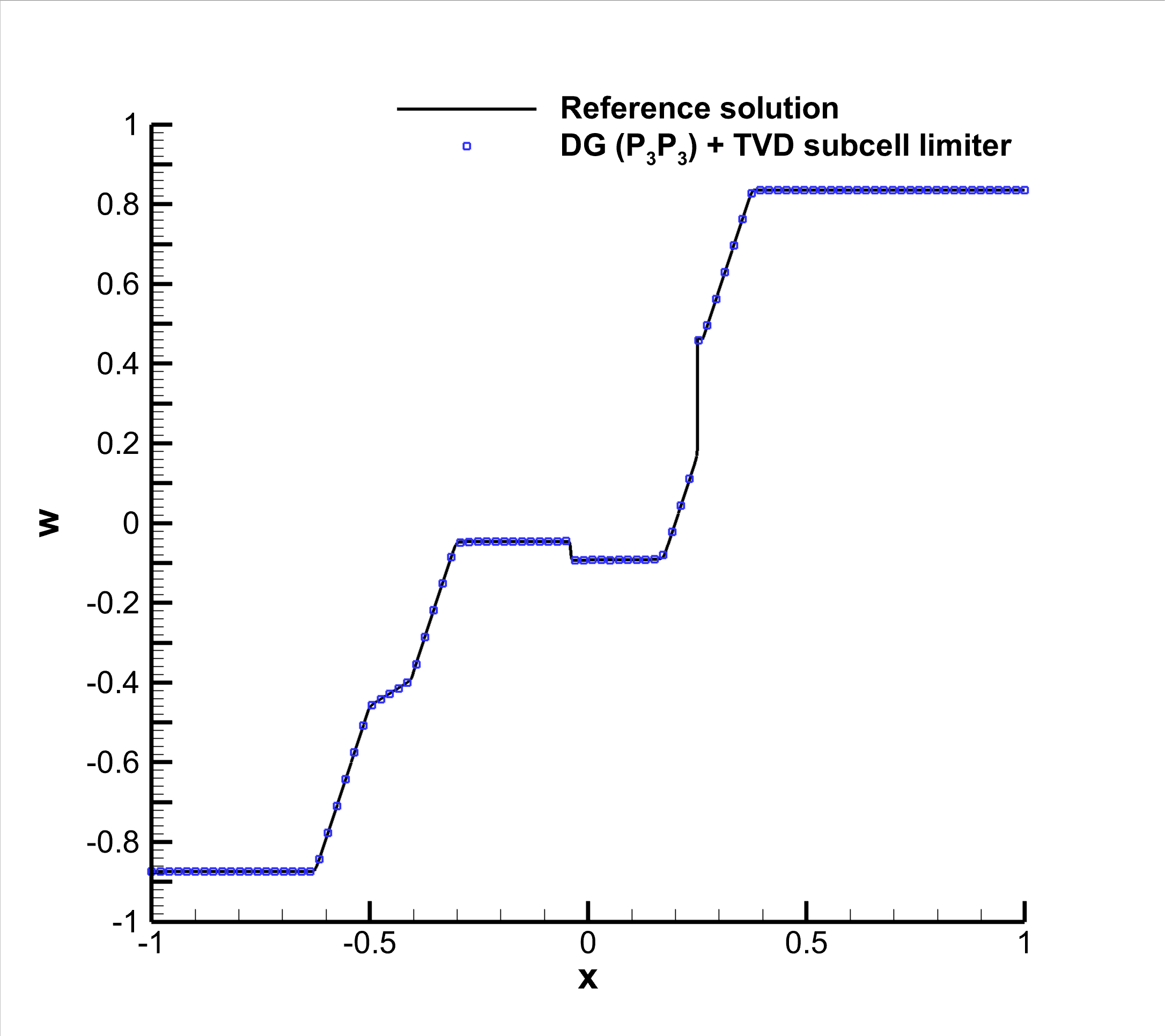}      
		\end{tabular}
		\caption{1D Riemann problem solved with a fourth-order ADER-DG scheme with TVD subcell limiter on a Cartesian mesh at time $t=0.25$. Top row: densities of each phase, $\rho_1$ and $\rho_2$. Second row: mixture density $\rho$ and $\alpha$. Third row: velocities $\vel_1$ and $\vel_2$. Bottom row: mixture velocity $\vel$ (left) and relative velocity $\urel=u_1-u_2$ (right).} 
		\label{fig:RP1D}
	\end{center}
\end{figure}

This section is devoted to studying the behavior of the proposed methodology in the presence of shocks. First, we solve one of the one-dimensional Riemann problems (RP) proposed in~\cite{Thein2022}, where a shock in one phase appears inside a rarefaction of the other phase. This RP presents a discontinuity in $x=0$ and has the left and right states shown in Table~\ref{tab:RP1}.
\begin{table}[!ht]
	\centering
    \caption{Left and right states of the RP1}
    \begin{tabular}{cccccccccc}
    \toprule
	 & $\alpha$ & $\rho_1$ & $\rho_2$ & $u_1$ & $u_2$  \\	\midrule 
$Q_L$ & 0.7 & 1.2449 & 1.2969 & -1.2638 & -0.38947 \\
$Q_R$ & 0.3 & 0.60312 & 0.73436 & 0.43059 & -0.40507 \\\bottomrule 
    \end{tabular}
    \label{tab:RP1}
\end{table} 
Since the problem that we want to reproduce is a 1D problem, $\vel^2=\vel^3=0$ and $\urel^2=\urel^3=0$. The computational domain is $\Omega = [-1,1]$ and has been discretized using a fourth-order ADER-DG scheme ($N=3$) with \aposteriori subcell limiter on a mesh with 1024 cells. The simulation is performed up to $t=0.25$, and the CFL number is set to $0.25$. Two ideal gases are considered for both phases with EOS~\eqref{eq:pressure_ideal}, setting $s_i=0$, $\gamma_1=1.4$, and $\gamma_2=2$, respectively. In Figure~\ref{fig:RP1D}, the numerical results are shown together with the reference solution computed with a second-order MUSCL-Hancock scheme based on the Rusanov flux as approximate Riemann solver and using a mesh spacing of $\Delta x = 0.5\cdot 10^{-4}$ (see~\cite{Thein2022} for further details). We observe an excellent agreement between our numerical solution and the reference solution. Looking at the density and velocity of each phase, we see that the rarefaction only affects the related phase, but that the interaction of the two rarefactions is observed in the mixing quantities.
In addition, in the density of the first phase, we can see that to the right of the contact a rarefaction begins (which does not affect the second phase) until the shock occurs. Once the shock appears, the density jumps according to the jump conditions. Then on the right a plateau of the right state of the shock is observed, and then the rarefaction continues again, see \cite{Thein2022} for a detailed discussion of this phenomenon.

\subsection{2D explosion problems} \label{sec:circular_explosion}

In this section, we solve the system for multi-phase flows in 2D in a circular computational domain with radius $R=1$. The initial condition is given by 
\begin{equation}
\bQ(\bx,t) = \left\{
			  \begin{array}{ll}
				\bQ_L & \mbox{if }|\bx|<0.5,\\
				\bQ_R & \mbox{otherwise },
			  \end{array}
			  \right.
			  \label{eq:CE}
\end{equation}
where $\bQ_L$ and $\bQ_R$ are described in Table~\ref{tab:CE}.
\begin{table}[!ht]
	\centering
	\caption{Left and right states of the circular explosion problems}
	\begin{tabular}{cccccccccc}
		 \toprule
		\multicolumn{10}{c}{CE1} \\ \midrule
		& $\alpha$ & $\rho_1$ & $\rho_2$ & $u^1_1$ & $u^2_1$ & $u^3_1$ & $u^1_2$ & $u^2_2$ & $u^2_2$ \\	\midrule 
		$\bQ_L$ & 0.4 & 2 & 1.5 & 0 & 0 & 0 & 0 & 0 & 0\\
		$\bQ_R$ & 0.8 & 1 & 0.5 & 0 & 0 & 0 & 0 & 0 & 0 \\\midrule
		 \midrule 
		\multicolumn{10}{c}{CE2} \\ \midrule
		& $\alpha$ & $\rho_1$ & $\rho_2$ & $u^1_1$ & $u^2_1$ & $u^3_1$ & $u^1_2$ & $u^2_2$ & $u^2_2$ \\	\midrule 
		$\bQ_L$ & 0.7 & 1 & 2 & 0 & 0 & 0 & 0 & 0 & 0\\
		$\bQ_R$ & 0.3 & 2 & 1 & 0 & 0 & 0 & 0 & 0 & 0 \\\bottomrule
	\end{tabular}
	\label{tab:CE}
\end{table} 
As a reference solution, we will solve the following equivalent (non-conservative) PDE in radial direction with geometric reaction source terms
\begin{subequations}\label{eq:barotropic_radial_CE}
\begin{align}
&\frac{\dalpha{1}}{\dt}+\velr \frac{\dalpha{1}}{\der r}=0,\label{eq:volume_frac_cyl_CE}\\
&\frac{\der\alpha_1\rho_1}{\dt}+\frac{\der (\alpha_1\rho_1 \velr_1)}{\der r}=-\frac{d}{r}(\alpha_1\rho_1 \velr_1),\label{eq:mass_density1_cyl_CE}\\
&\frac{\der\alpha_2\rho_2}{\dt}+\frac{\der (\alpha_2\rho_2 \velr_2)}{\der r}=-\frac{d}{r}(\alpha_2\rho_2 \velr_2),\label{eq:mass_density2_cyl_CE}\\
&\frac{\der\rho \velr}{\dt}+ \frac{\der(\alpha_1\rho_1 (\velr_1)^2+ \alpha_2\rho_2 (\velr_2)^2 + p )}{\der r}=-\frac{\alpha_1\rho_1 (u_1^r)^2}{r}-\frac{\alpha_2\rho_2 (u_2^r)^2}{r},\label{eq:momentum_cyl_CE}\\
&\frac{\der \urelr}{\dt}+\frac{\der}{\der r} \left( \frac{1}{2}(\velr_1)^2-\frac{1}{2}(\velr_2)^2 + h_1-h_2\right)=0,\label{eq:relative_velocity_cyl_CE}
\end{align}
\end{subequations}
where the parameter $d$ is the number of spatial dimensions minus one.

The 2D computations have been performed using a fourth-order ($N=3$) ADER-DG scheme with \aposteriori subcell limiter. The computational domain is $\Omega=[-1,1]\times [-1,1]$ and has been discretized using a Cartesian mesh with $256\times 256$ elements.  Following~\eqref{eq:CE}, the left state of the RP has been taken as the inner state and the right state of the same RP as the outer state. The reference solution has been computed by solving~\eqref{eq:barotropic_radial_CE} with 128000 cells using a second-order TVD finite volume method with the Rusanov flux. The simulation is performed up to $t=0.1$ with two ideal gases, so for the two phases, the EOS is given by~\eqref{eq:pressure_ideal}, with $s_i=0$, $\gamma_1=1.4$, and $\gamma_2=2$, respectively. 

\begin{figure}
	\begin{center}
		\begin{tabular}{cc} 
			\includegraphics[trim= 5 5 5 5,clip,width=0.35\textwidth]{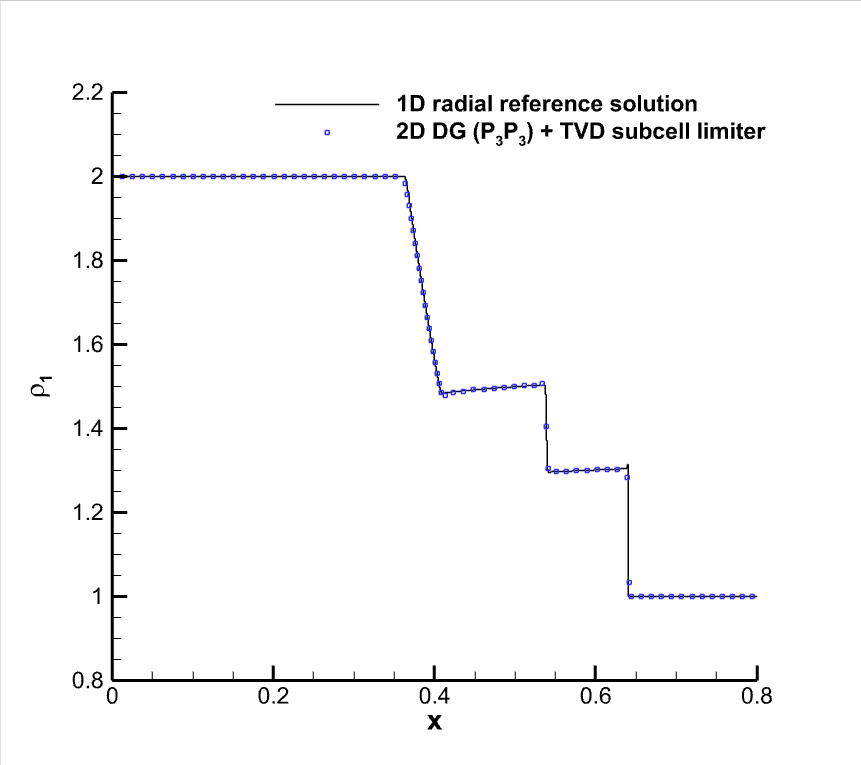}   & 
			\includegraphics[trim= 5 5 5 5,clip,width=0.35\textwidth]{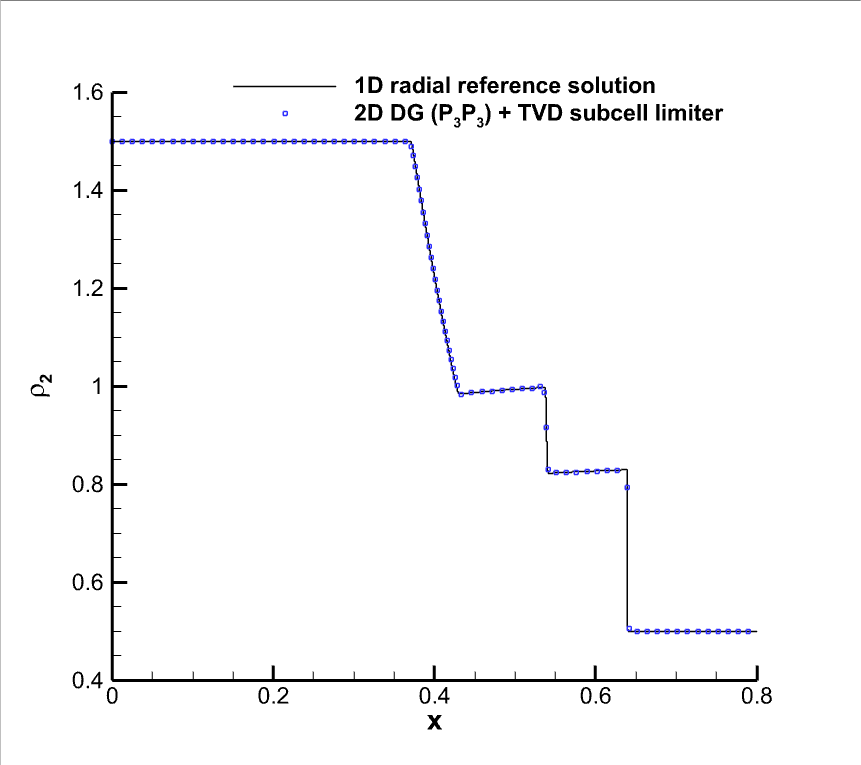}   \\[-0.2cm] 
			\includegraphics[trim= 5 5 5 5,clip,width=0.35\textwidth]{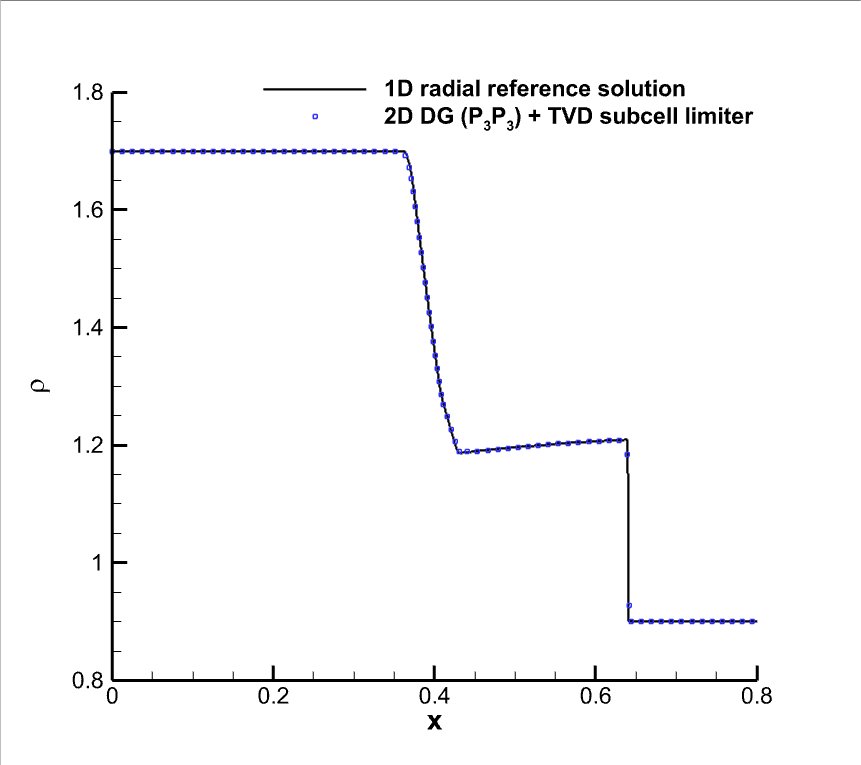}    &   
			\includegraphics[trim= 5 5 5 5,clip,width=0.35\textwidth]{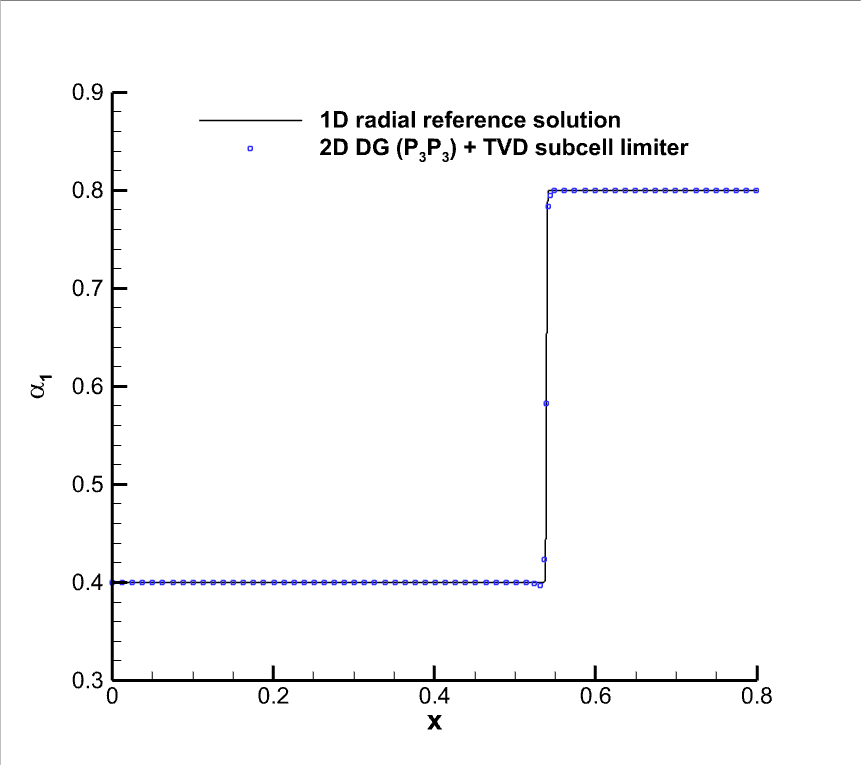}  \\[-0.2cm]
			\includegraphics[trim= 5 5 5 5,clip,width=0.35\textwidth]{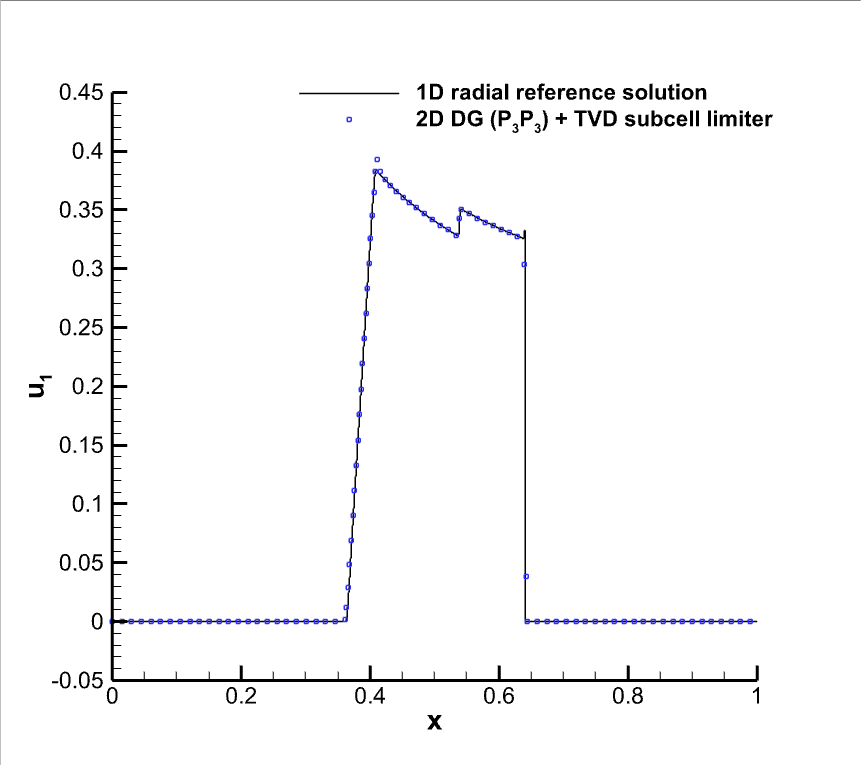}     & 
			\includegraphics[trim= 5 5 5 5,clip,width=0.35\textwidth]{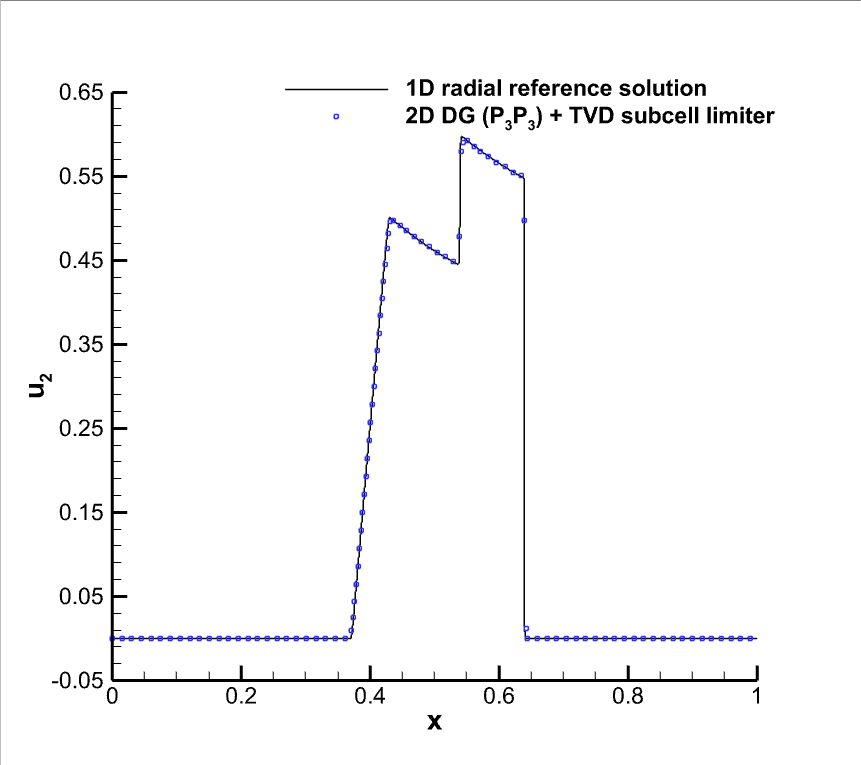}     \\[-0.2cm]
			\includegraphics[trim= 5 5 5 5,clip,width=0.35\textwidth]{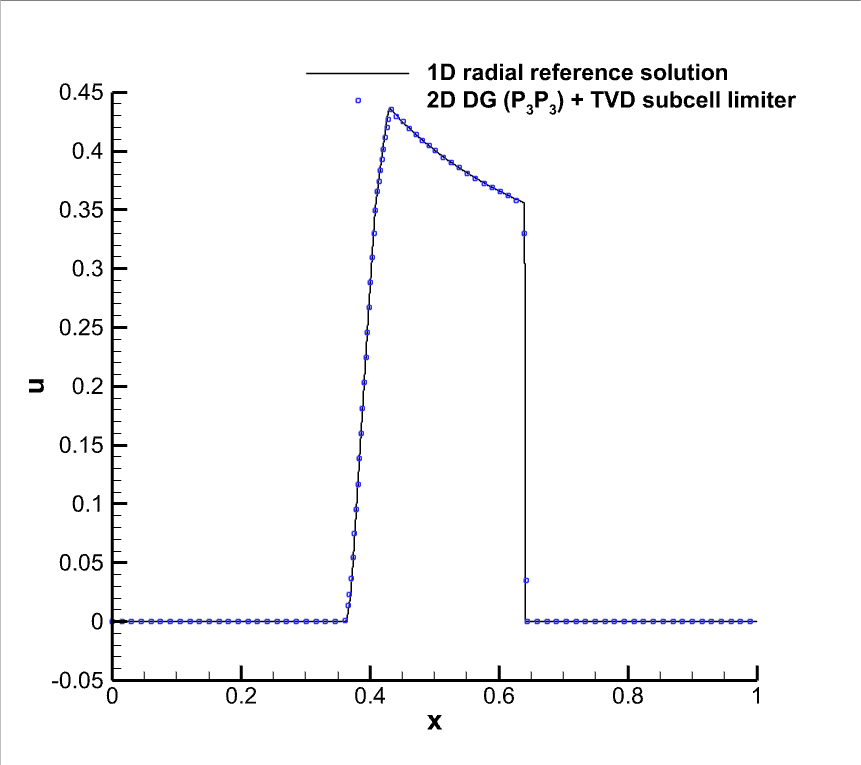}      &
			\includegraphics[trim= 5 5 5 5,clip,width=0.35\textwidth]{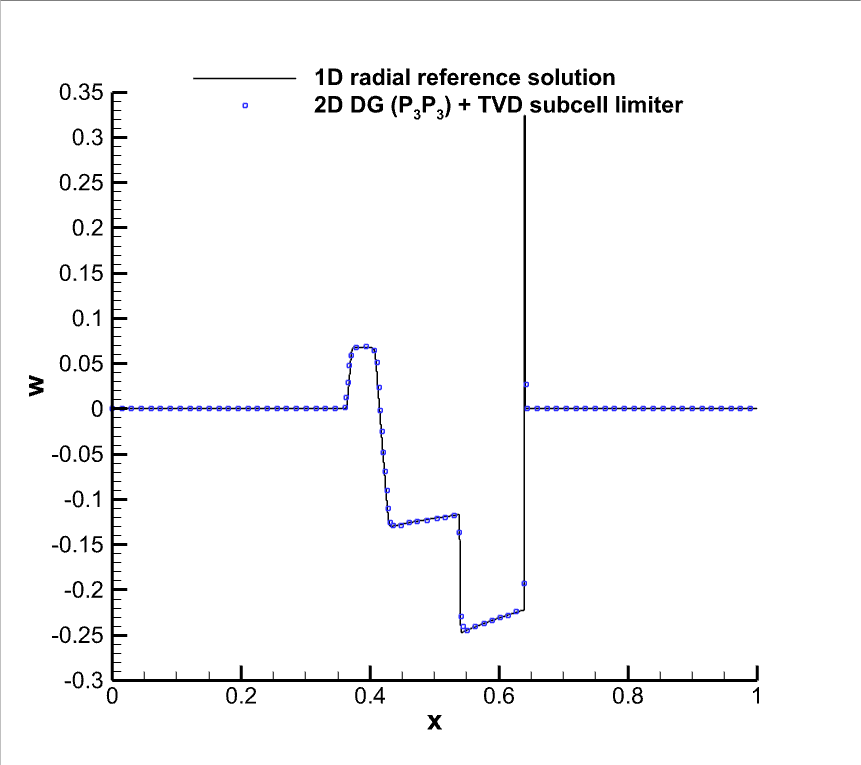}      
		\end{tabular} 
		\caption{2D circular explosion problem for initial condition CE1 in~\ref{tab:CE} solved on a Cartesian mesh at time $t=0.1$, in comparison with the radial reference solution. Top row: densities of each phase, $\rho_1$ and $\rho_2$. Second row: mixture density $\rho$ and $\alpha$. Third row: velocities $\vel_1$ and $\vel_2$. Bottom row: mixture velocity $\vel$ (left) and relative velocity $\urel=u_1-u_2$ (right).} 
		\label{fig:CE1}
	\end{center}
\end{figure}

\begin{figure}
	\begin{center}
		\begin{tabular}{cc} 
			\includegraphics[trim= 5 5 5 5,clip,width=0.35\textwidth]{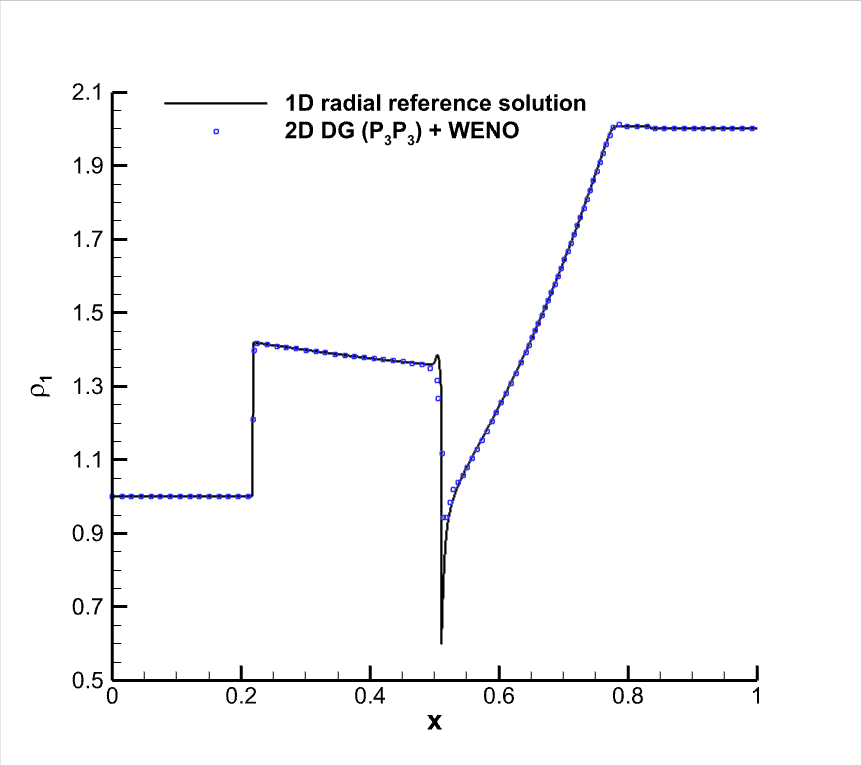}   & 
			\includegraphics[trim= 5 5 5 5,clip,width=0.35\textwidth]{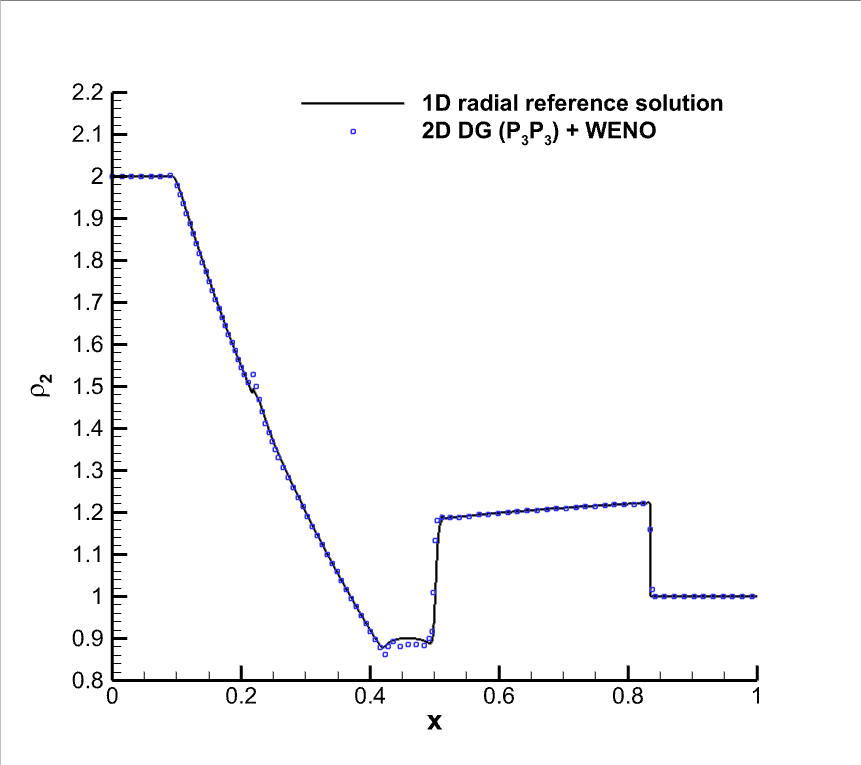}   \\[-0.2cm] 
			\includegraphics[trim= 5 5 5 5,clip,width=0.35\textwidth]{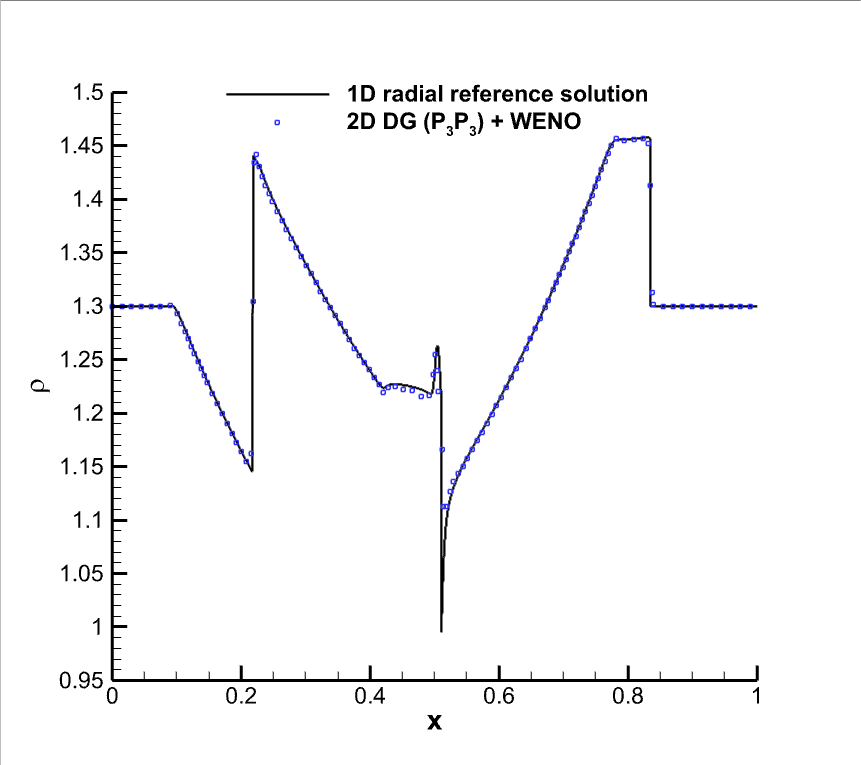}    &   
			\includegraphics[trim= 5 5 5 5,clip,width=0.35\textwidth]{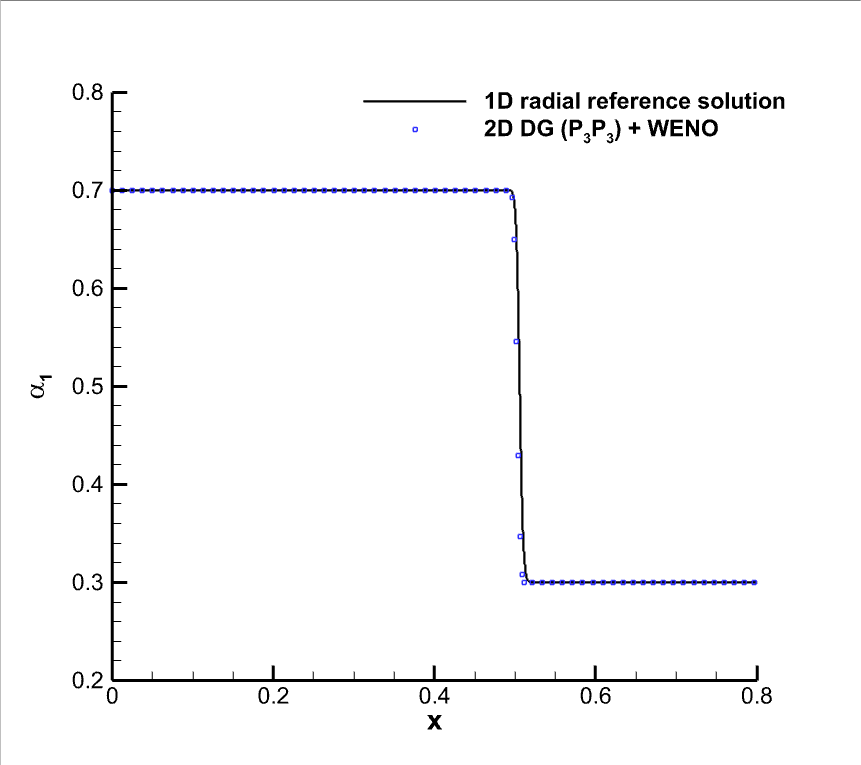}  \\[-0.2cm]
			\includegraphics[trim= 5 5 5 5,clip,width=0.35\textwidth]{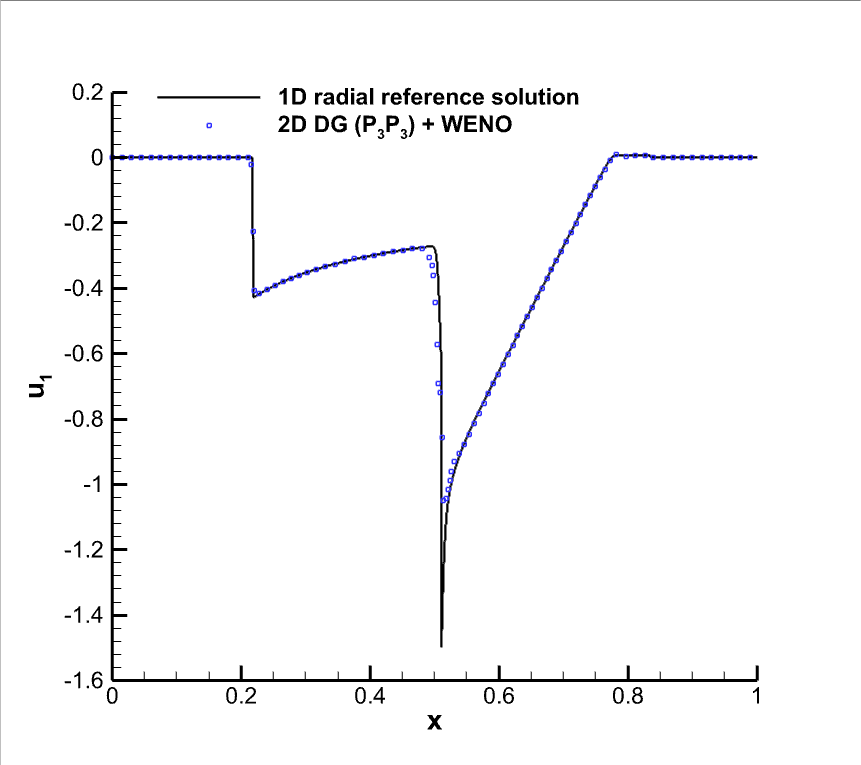}     & 
			\includegraphics[trim= 5 5 5 5,clip,width=0.35\textwidth]{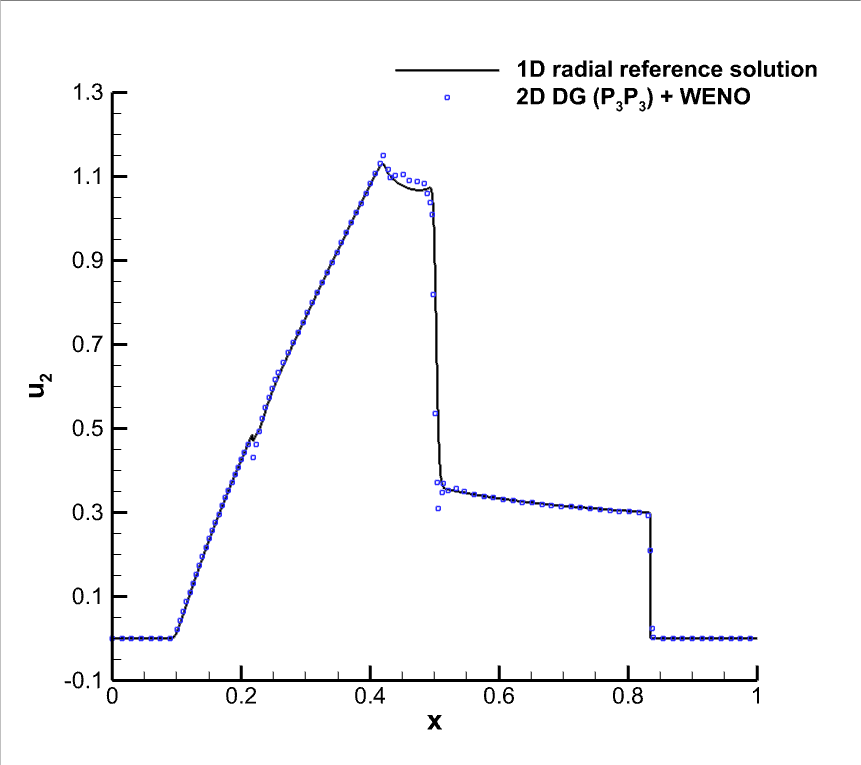}     \\[-0.2cm]
			\includegraphics[trim= 5 5 5 5,clip,width=0.35\textwidth]{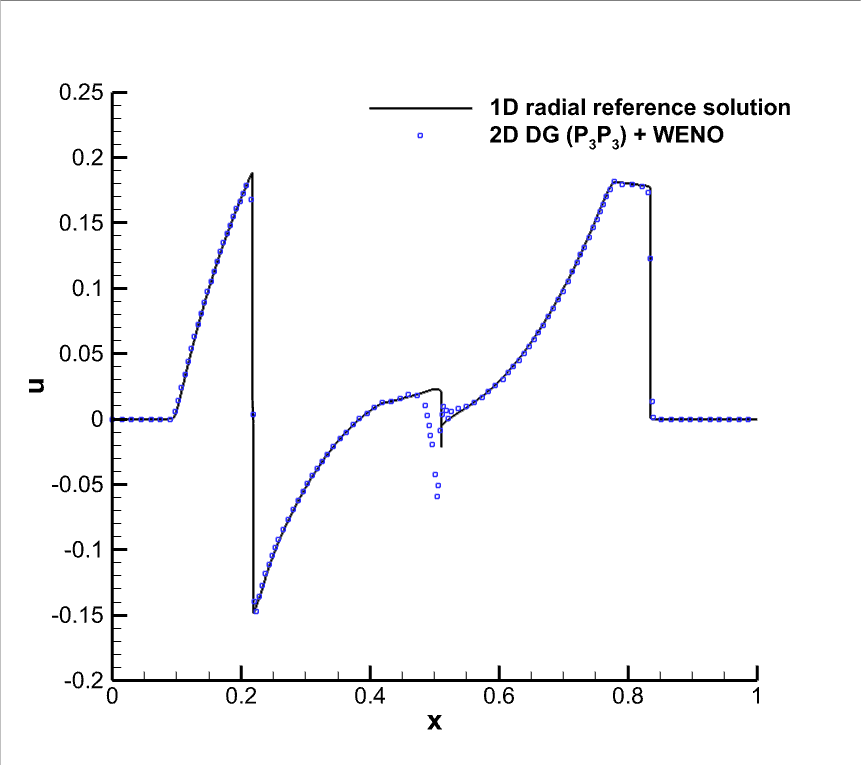}      &
			\includegraphics[trim= 5 5 5 5,clip,width=0.35\textwidth]{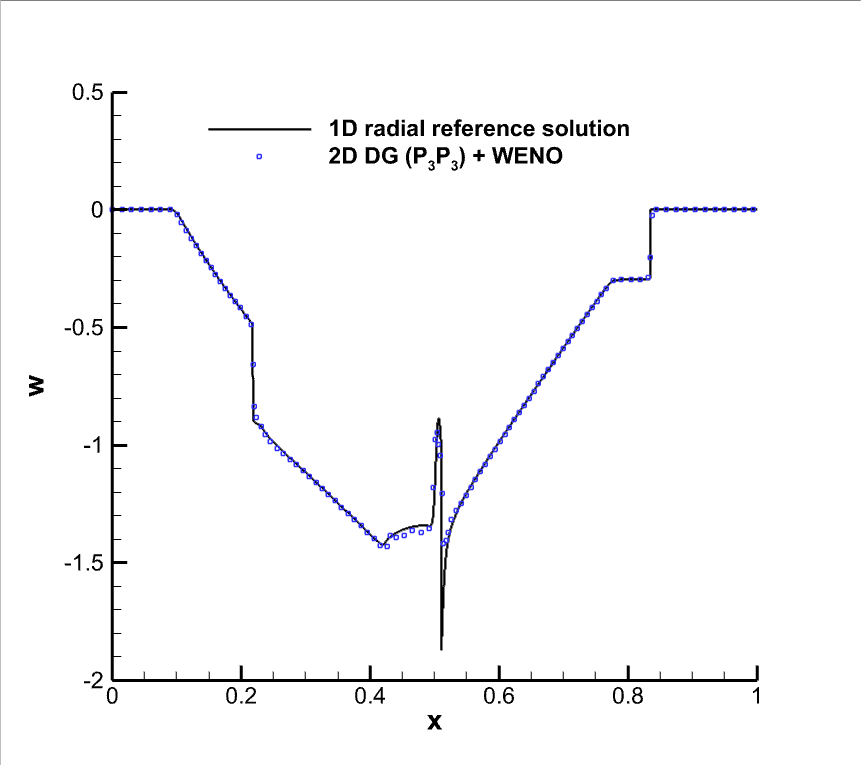}      
		\end{tabular}
		\caption{2D circular explosion problem for initial condition CE2 in~\ref{tab:CE} solved on a Cartesian mesh at time $t=0.2$, compared with the radial reference solution. Top row: densities of each phase, $\rho_1$ and $\rho_2$. Second row: mixture density $\rho$ and $\alpha$. Third row: velocities $\vel_1$ and $\vel_2$. Bottom row: mixture velocity $\vel$ (left) and relative velocity $\urel=u_1-u_2$ (right).} 
		\label{fig:CE2}
	\end{center}
\end{figure}

\begin{figure}
	\begin{center}
		\begin{tabular}{cc} 
			\includegraphics[trim= 5 5 5 5,clip,width=0.45\textwidth]{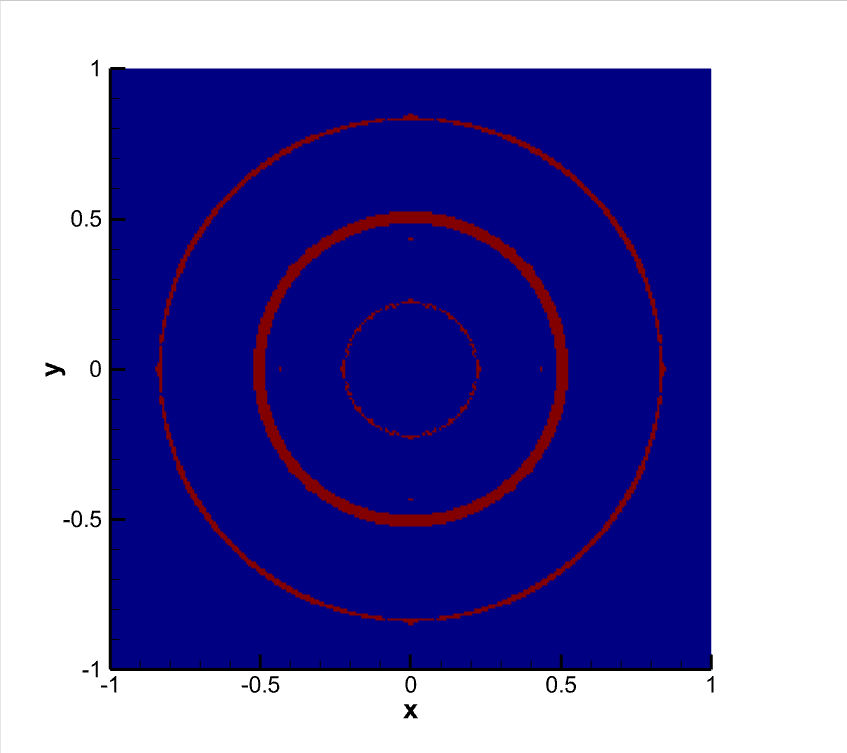}   & 
			\includegraphics[trim= 5 5 5 5,clip,width=0.45\textwidth]{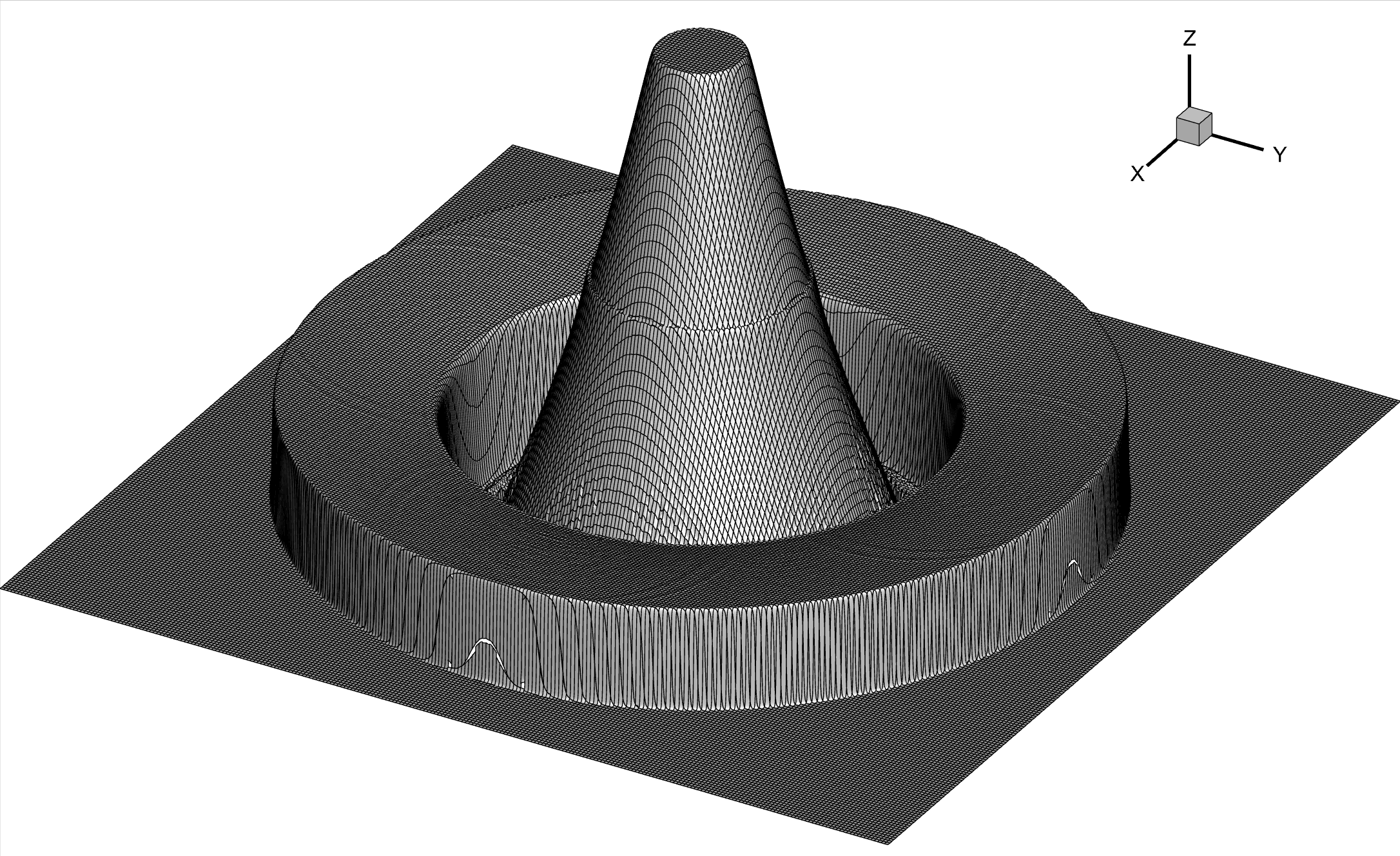} 
		\end{tabular}
		\caption{Left: Limiter map of the explosion problem in 2D. The values in red mean that the limiter is activated. Right: 3D plot with the variable $\rho_2$ in the $z-$axis.} 
		\label{fig:CE2_limiter map}
	\end{center}
\end{figure}

Figures~\ref{fig:CE1} and~\ref{fig:CE2} show the numerical results of two circular explosion problems, with the initial conditions of Table~\ref{tab:CE} and a final time of $t=0.1$ for the first one and $t=0.2$ for the second one. The numerical solution obtained with the ADER-DG method is then compared with the radial reference solution, showing excellent agreement.  
Moreover, Figure~\ref{fig:CE2_limiter map} shows the limiter map of the second explosion problem. The values highlighted in blue are those DG elements where the limiter is not activated, and the red ones are the troubled zones where the \aposteriori subcell FV limiter is activated.

\subsection{Dambreak problem} \label{sec:dambreak} 
Finally, a two-dimensional dambreak problem is solved using the barotropic two-phase model studied in this paper. In this case, the source term included in the momentum equation~\eqref{eq:cons_momentum} is non-zero, as a non-zero gravity source is considered, hence $\bg=(0,-g,0)$ with $g=9.81$. The computational domain is $\Omega=[0,4]\times[0,2]$, where the water domain is $\Omega_2=[0,2]\times[0,1]$ and the air domain is given by $\Omega_1=\Omega \backslash \Omega_2$. 

\begin{figure}
	\begin{center} 
	\includegraphics[trim=2 0 2 2,clip,width=0.75\textwidth]{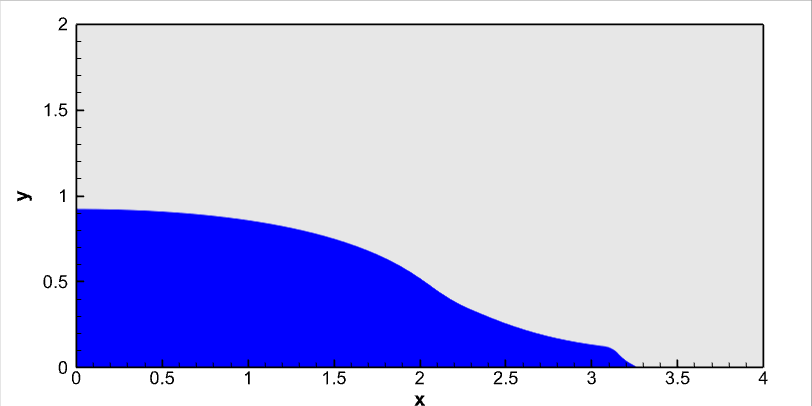}  \\ 
	\includegraphics[trim=2 0 2 2,clip,width=0.75\textwidth]{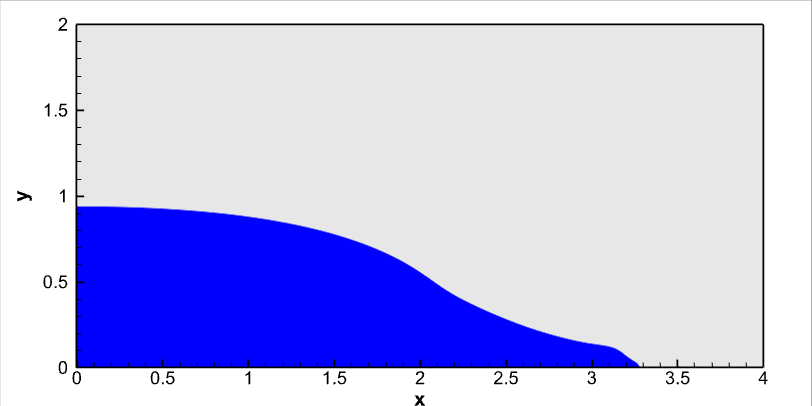}   \\ 
	\includegraphics[trim=2 2 2 2,clip,width=0.75\textwidth]{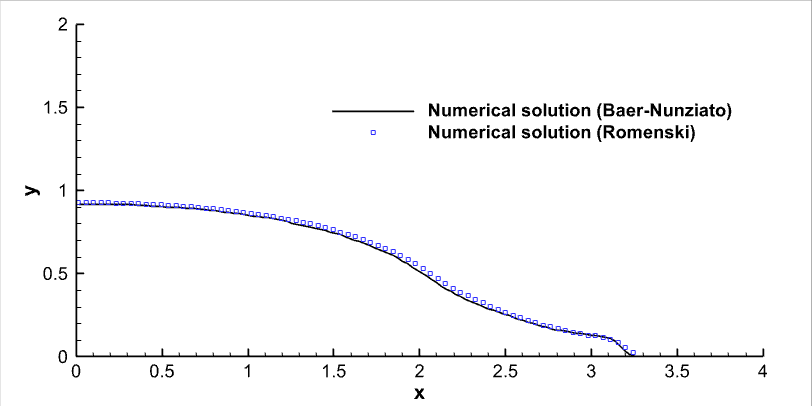}
		\caption{Dambreak problem at time $t = 0.4$. Top: reference solution, computed with a third-order ADER-WENO finite volume scheme on a very fine uniform Cartesian grid, solving the inviscid and barotropic reduced Baer--Nunziato model presented in~\cite{Dumbser2011TwoPhase}. Center: Numerical solution, computed using an ADER-DG scheme with \aposteriori sub-cell limiter, to solve the barotropic SHTC model proposed in this work. Bottom: Comparison of the free-surface profile obtained for both models.}
		\label{fig:dambreak}
	\end{center}
\end{figure}

The domain $\Omega$ has been discretized with a uniform Cartesian mesh with 256x128 cells, using an ADER-DG scheme with $N=3$ and \aposteriori subcell FV limiter. The simulation has been performed until a final time of $t=0.4$, and a slip wall boundary condition is imposed on all boundaries. Following Section~\ref{sec:eos}, an ideal gas is considered in $\Omega_1$, i.e., the EOS is given by~\eqref{eq:energy_ideal}, with parameters $c_{01}=1$, $\gamma_1 = 1.4$, $\rho_{01}=1$, $\alpha=\varepsilon$. The initial pressure profile is assumed hydrostatic, $p=\rho_{01} g (y-2)$. The EOS for the liquid is a stiffened gas EOS given by~\eqref{eq:energy_stiffened} where $c_{02}=20$, $\gamma_2 = 2$, $\rho_{02}=1000$, $\alpha=1-\varepsilon$, and again a hydrostatic pressure profile $p=\rho_{02} g (y-1)$ is imposed initially. The simulation was performed with $\varepsilon=0$, i.e., initially the phase volume fractions are really set to zero and unity, respectively. To obtain the value of the primitive variable $\rho_k$ it is necessary to divide by $\alpha_k$, and in this simulation, there exist areas with $\alpha_k=0$, and it is necessary to apply a \textit{filter} that avoids division by zero. In this paper, the density variables are filtered as follows, 
\begin{equation*}
\rho_k   = \frac{\rho_k \alpha_k^2 + \rho_{0k} \epsilon }{\alpha_k^2 + \epsilon},  
\end{equation*}
see also \cite{Tavelli2018}, and the filter parameter is set $\epsilon=10^{-12}$. 
The numerical results have been compared with the solution of the reduced barotropic Baer--Nunziato model given in \cite{Dumbser2011TwoPhase}. Figure~\ref{fig:dambreak} shows the values obtained at time $t=0.4$ calculated with the reduced Baer--Nunziato model in the upper plot, the solution calculated with the method proposed in this paper in the center, and a direct comparison between both models in the bottom plot, showing an excellent agreement between both models. Similar results have also been recently obtained with a novel Arbitrary--Lagrangian--Eulerian hybrid finite volume / finite element method applied to the incompressible Navier--Stokes equations on moving unstructured meshes, see~\cite{Busto2023HybridALE}.  

\section{Conclusion}
\label{sec:conclusion} 
In this paper, the barotropic version of the conservative SHTC model for compressible two-fluid flows of Romenski \textit{et al.} has been solved for the first time using high-order ADER discontinuous Galerkin schemes in combination with an \aposteriori subcell FV limiter. Since the model is only weakly hyperbolic in the general multidimensional case, two different methodologies have been presented to restore the strong hyperbolicity: i) a generalized Lagrangian multiplier (GLM) curl-cleaning approach and ii) the addition of the Godunov-Powell terms to symmetrize the system. We obtain a full set of linearly independent eigenvectors with both methodologies, proving that strong hyperbolicity can indeed be restored.  

A high-order ADER discontinuous Galerkin finite element scheme with \aposteriori subcell finite volume limiter has been used to deal with discontinuities and steep gradients in the solutions. To validate the model and the proposed method, a numerical convergence analysis has been carried out, and the high order of the method has been confirmed. For this purpose, we have constructed a new exact analytical and stationary equilibrium solution of the PDE system in cylindrical coordinates. Then, several Riemann problems in one and two dimensions have been simulated to show the behavior of the proposed methodology in the presence of shocks. First, a 1D Riemann problem where a shock in one phase appears inside the rarefaction of the other phase has been simulated. The results have been compared with those presented in~\cite{Thein2022}, showing an excellent agreement. Then two 2D explosion problems were solved. Thanks to the radial symmetry of the problem, the obtained results have been compared with an equivalent 1D reference solution, showing the accuracy of the proposed methodology even in presence of sharp gradients in the solution. Finally, a dambreak test case has been considered, where the initial values of the volume fractions are set to $\alpha_1=0$ and $\alpha_2=1$. The numerical results are compared with those obtained for a reduced barotropic Baer--Nunziato-type model, showing an excellent agreement between both models. 

As future work, we plan to extend our methodology to compressible multi-phase flows with more than two phases and, in addition, to include also solids governed by the equations of nonlinear hyperelasticity, see, e.g.~\cite{GPRmodel,Romenski2020,Romenski2022}. Furthermore, we will also apply exactly curl-free methods to the two-phase model discussed in this paper, such as the curl-free schemes recently forwarded in \cite{Balsara2023CurlFree,Boscheri2021SIGPR,Dhaouadi2023expspnsk}.

\section*{Acknowledgments}
This research was funded by the Italian Ministry of Education, University and Research (MIUR) in the frame of the Departments of Excellence  Initiative 2018--2027 attributed to DICAM of the University of Trento (grant L. 232/2016) and in the frame of the PRIN 2017 project \textit{Innovative numerical methods for evolutionary partial differential equations and  applications}.
L.R. acknowledges funding from the Spanish Ministry of Universities and the European Union-Next GenerationEU under the project RSU.UDC.MS15. MD and LR are members of the GNCS group of INdAM.
The authors would like to acknowledge support from the CESGA, Spain, for the access to the FT3 supercomputer and to the CINECA award under the ISCRA initiative, for the availability of high performance computing resources and support (project number IsCa3\_NuMFluS). 
L.R. gratefully acknowledges Dr. Firas Dhaouadi and Dr. Ilya Peshkov for the interesting discussions and support that have allowed the successful development of this work.

\section*{Dedication}

This paper is dedicated to Gerald Warnecke at the occasion of his 65$^{\textnormal{th}}$ birthday and in honor of his groundbreaking scientific contributions to the field of numerical methods for hyperbolic PDE. 
The authors are also very grateful for the friendship and all the inspiring discussions over the years.

\section*{Conflict of Interest}

The authors declare that they have no conflict of interest. 

\section*{Data availability}

The data can be obtained from the authors on reasonable request. 

\bibliographystyle{elsarticle-num}
\bibliography{./mibiblio}

\end{document}